\documentclass[a4paper,11pt]{amsproc}
\usepackage[margin=3.6cm]{geometry}
\usepackage{amsmath,amsfonts,enumerate,amssymb,amsthm,color}
\usepackage[dvips]{epsfig}
\usepackage{pgf,tikz}
\usepackage{fancyhdr}

\title{Automorphism groups of rational circulant graphs through the use of Schur rings}
\author{Mikhail Klin \and Istv\'an Kov\'acs}
\address{Department of Mathematics \& Computer Science, Ben-Gurion University of the Negev, 
84105 Beer-Sheva, Israel}
\email{klin@cs.bgu.ac.il}
\address{FAMNIT, University of Primorska, Glagolja\v{s}ka 8, 6000 Koper,  Slovenia}
\email{istvan.kovacs@upr.si}
\thanks{Second 
author was supported in part by ARRS -- 
Agencija za raziskovanje Republike Slovenije, program no. P1-0285.}

\newtheorem{lem}{Lemma}[section]
\newtheorem{prop}[lem]{Proposition}
\newtheorem{cor}[lem]{Corollary}
\newtheorem{thm}[lem]{Theorem}

\theoremstyle{definition}

\newtheorem{exm}[lem]{Example}
\newtheorem{defi}[lem]{Definition}

\def\N{\mathbb{N}}
\def\Q{\mathbb{Q}}

\def\Z{\mathbb{Z}}
\def\A{\mathcal{A}}
\def\B{\mathcal{B}}
\def\X{\mathfrak{X}}
\def\F{\mathcal{F}}

\def\P{\mathcal{P}}
\def\un{\underline}
\def\ov{\overline}
\def\trS{\stackrel{\; \; \circ}{S}}

\def\trA{\stackrel{\; \; \; \circ}{\A}}

\def\proof{\noindent{\sc Proof}.\ }

\def\QED{\hfill\vrule height .9ex width .8ex depth -.1ex \bigskip}

\DeclareMathOperator{\cay}{Cay}
\DeclareMathOperator{\aut}{Aut}

\DeclareMathOperator{\sym}{Sym}

\DeclareMathOperator{\co}{As}
\DeclareMathOperator{\orb}{2-Orb}
\DeclareMathOperator{\anc}{Anc}
\DeclareMathOperator{\pbs}{\mathsf{PBS}}

\DeclareMathOperator{\basic}{\mathsf{Basic}}

\newcommand{\sg}[1]{\langle {#1}\rangle}
\newcommand{\sgg}[1]{\langle\!\langle {#1}  \rangle\!\rangle}
\newcommand{\comment}[1]{}
\newcommand{\subl}[1]{_{[#1]}}
\newcommand{\supl}[1]{^{[#1]}}

\begin{document}

\begin{abstract}
The paper concerns the automorphism groups of Cayley graphs over cyclic groups which have a
rational spectrum (rational circulant graphs for short). With the aid of the techniques of Schur rings it is shown that the problem is equivalent
to consider the automorphism groups of orthogonal group block structures of cyclic groups.
Using this observation, the required groups are expressed in terms of generalized
wreath products of symmetric groups. 
\end{abstract}

\maketitle

\pagestyle{fancy}
\fancyhead{}
\fancyhead[CE]{\sc mikhail klin and istv\'{a}n kov\'{a}cs}
\fancyhead[CO]{\sc automorphism groups of rational circulant graphs}
\fancyfoot{}
\fancyhead[RO,LE]{\thepage}

\section{Introduction}

A circulant graph with $n$ vertices is a Cayley graph over the cyclic group $\Z_n$, 
i.e., a graph having an automorphism which permutes all the vertices into a full cycle. There is a vast literature investigating
various properties of this class of graphs. In this paper we focus on their automorphisms. By the definition,
the automorphism groups contain a regular cyclic subgroup. The study of permutation groups with a regular cyclic subgroup goes back to the work of Burnside and Schur. Schur proved that if the group is primitive of
composite degree, then it is doubly transitive (see \cite{Sch33}). The complete list of such primitive groups was 
given recently by the use of the classification of finite simple groups, see \cite{Jon02,Li03}.

One might expect transparent descriptions of the automorphism groups of circulant graphs by restricting to a 
suitably chosen family. A natural restriction can be done with respect to the order $n$ of the graph. For instance, we 
refer to the papers \cite{DobM05,KliP,Kov} dealing with the case when $n$ is a square-free number, $n=p^e$ ($p$ is an odd prime), and $n=2^e$, respectively.
In the present paper we choose another natural family by requiring the graphs to have a rational
spectrum, i.e., the family of rational circulant graphs.


To formulate our main result some notations are in order. For $n\in \N$, we let $[n]$ denote
the set $\{1,\ldots,n\}$, and $S_n$ the group of all permutations of $[n]$. Let
$([r],\preceq)$ be a poset on $[r]$. We say that $([r],\preceq)$ is
{\em increasing} if $i \preceq j$ implies $i \le j$ for all $i,j \in [r]$. Below
$\prod_{([r],\preceq)}S_{n_i}$ denotes the generalized wreath
product, defined by $([r],\preceq)$ and the groups $S_{n_1},\ldots,S_{n_r}$,
acting on the set $[n_1] \times \cdots \times [n_r]$. For
the precise formulation, see Definition~\ref{def:GWP}.

\medskip

Our main result is the following theorem.

\begin{thm}\label{thm:main}
Let $G$ be a permutation group acting on the cyclic group $\Z_n$, $n \ge 2$.
The following are equivalent:
\begin{enumerate}[(i)]
\item $G=\aut(\cay(\Z_n,S))$ for some rational circulant graph $\cay(\Z_n,S)$.
\item $G$ is a permutation group, which is permutation isomorphic to a generalized wreath product
$\prod_{([r],\preceq)}S_{n_i}$, where $([r],\preceq)$ is an increasing poset, and
$n_1,\ldots,n_r$ are in $\N$ satisfying
\end{enumerate}

\begin{enumerate}[\qquad (a)]
\item $n=n_1 \cdots n_r$,
\item $n_i \ge 2$ for all $i \in \{1,\ldots,r\}$,
\item $(n_i,n_j)=1$ for all $i,j \in \{1,\ldots,r\}$ with $i \not\preceq j$.
\end{enumerate}
\end{thm}

To the number $n_i$ in (ii) we shall also refer to as the {\em weight} of 
node $i$ in the poset $([r],\preceq)$.
The following examples serve as illustrations of Theorem~\ref{thm:main}.

\begin{exm}\label{ex:n=6}
Here $n=6$. Up to complement, there are four rational circulant graphs:
$$ K_6, \; K_2 \times K_3, \; K_{3,3}, \; K_{2,2,2}.$$
The corresponding automorphism groups: $S_6, S_2 \times S_3, S_2 \wr S_3$, and
$S_3 \wr S_2$.

In part (ii) we get $G=S_6$ for $r=1$. If $r=2$, then any choice $n_1,n_2\in \{2,3\}$
with $n_1n_2=6$ gives weights of an increasing poset on $\{1,2\}$.
For instance, if $n_1=2$, $n_2=3$, and $([2],\preceq)$ is an anti-chain, then
$G=S_2 \times S_3$, and the same group is obtained if we switch the
values of weights. Changing the poset $([2],\preceq)$ to a chain we get the wreath products
$S_2 \wr S_3$ and $S_3 \wr S_2$. \QED
\end{exm}

\begin{exm}\label{exm:n=12}
Here $n=12$. In this example we consider the groups that can be derived from part (ii).
We have $G=S_{12}$ if $r=1$. If $r=2$, then similarly to the previous
example we deduce that $G$ is one of the groups:
$S_3 \times S_4$, $S_a \wr S_{n/a}$, $a \in \{2,3,4,6\}$.

Let $r=3$. The three nodes of $([3],\preceq)$ get weights $2,2,3$ by (a)-(b), and
because of (c) the two nodes with weight $2$ must be related.
The possible increasing posets are depicted in Figure~1.

\begin{figure}[h!]
\centering
\begin{tikzpicture}[scale=1.4] ---
\draw (0,0) circle (3pt);
\draw (0,0) node {\footnotesize $1$};
\draw (0.5,0) circle (3pt);
\draw (0.5,0) node {\footnotesize $2$};
\draw (0.5,1) circle (3pt);
\draw (0.5,1) node {\footnotesize $3$};
\draw (0.5,0.1) -- (0.5,0.9);
\draw (0.25,-1) node {\textrm (i)};

\draw (1.5,0) circle (3pt);
\draw (1.5,0) node {\footnotesize $2$};
\draw (2,0) circle (3pt);
\draw (2,0) node {\footnotesize $1$};
\draw (2,1) circle (3pt);
\draw (2,1) node {\footnotesize $3$};
\draw (2,0.1) -- (2,0.9);
\draw (1.75,-1) node {\textrm (ii)};

\draw (3,0) circle (3pt);
\draw (3,0) node {\footnotesize $3$};
\draw (3.5,0) circle (3pt);
\draw (3.5,0) node {\footnotesize $1$};
\draw (3.5,1) circle (3pt);
\draw (3.5,1) node {\footnotesize $2$};
\draw (3.5,0.1) -- (3.5,0.9);
\draw (3.25,-1) node {\textrm (iii)};

\draw (4.5,0) circle (3pt);
\draw (4.5,0) node {\footnotesize $1$};
\draw (5,0) circle (3pt);
\draw (5,0) node {\footnotesize $2$};
\draw (4.75,1) circle (3pt);
\draw (4.75,1) node {\footnotesize $3$};
\draw (4.5,0.1) -- (4.75,0.9);
\draw (5,0.1) -- (4.75,0.9);
\draw (4.75,-1) node {\textrm (iv)};

\draw (6,1) circle (3pt);
\draw (6,1) node {\footnotesize $2$};
\draw (6.25,0) circle (3pt);
\draw (6.25,0) node {\footnotesize $1$};
\draw (6.5,1) circle (3pt);
\draw (6.5,1) node {\footnotesize $3$};
\draw (6.25,0.1) -- (6,0.9);
\draw (6.25,0.1) -- (6.5,0.9);
\draw (6.25,-1) node {\textrm (v)};

\draw (7.5,0) circle (3pt);
\draw (7.5,0) node {\footnotesize $1$};
\draw (7.5,1) circle (3pt);
\draw (7.5,1) node {\footnotesize $2$};
\draw (7.5,2) circle (3pt);
\draw (7.5,2) node {\footnotesize $3$};
\draw (7.5,0.1) -- (7.5,0.9);
\draw (7.5,1.1) -- (7.5,1.9);
\draw (7.5,-1) node {\textrm (vi)};

\end{tikzpicture}
\caption{Increasing posets on $\{1,2,3\}$.}
\end{figure}
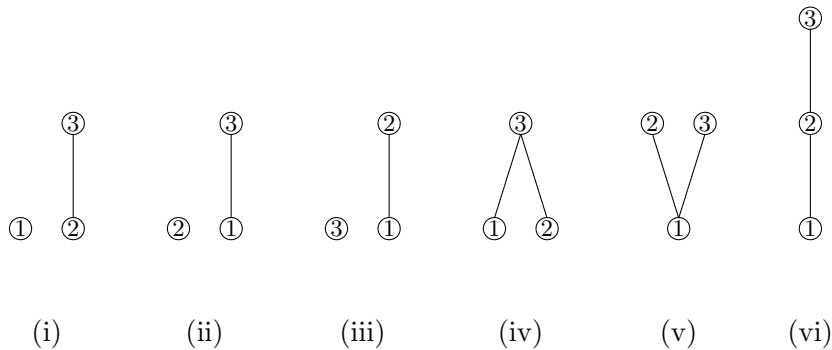

\noindent The weights are unique for posets (i)-(iii). 
In poset (iv) the only restriction is that $n_3=2$, in poset (v) 
the only restriction is that $n_1=2$, and weights
are arbitrarily distributed for poset (vi).
By Definition~\ref{def:GWP}, we obtain the following groups:

\begin{itemize} 
\item $S_3 \times (S_2 \wr S_2)$ corresponding to posets (i)-(iii),
\item $S_{2} \wr (S_2 \times S_3)$ corresponding to poset (iv)
\item $(S_2 \times S_3) \wr S_2$ corresponding to poset (v),
\item $S_3 \wr S_2 \wr S_2$, $S_2 \wr S_3 \wr S_2$ and
$S_2 \wr S_2 \wr S_3$ corresponding to poset (vi) 
(here the group depends also on the weights).
\end{itemize}

Finally, altogether we obtain exactly $12$ possible distinct groups (including the largest $S_{12}$ and 
the smallest of order $48$). Each such group appears exactly ones. (Attribution of the same group of order 
$48$ to three posets is an artifice, which results from the way of the presentation.)
Observe that, each of these groups is obtained using iteratively direct or wreath product
of symmetric groups. \QED
\end{exm}

For larger values of $n$ it is not true that generalized wreath product of symmetric groups may be obtained by an iterative use of direct and wreath products of symmetric groups.  An example of such situation appears for $n=36$,
and it will be discussed later on in the text. 

In deriving Theorem~\ref{thm:main} we follow an approach suggested by
Klin and P\"oshcel in \cite{KliP81}, which is to explore the Galois correspondence 
between overgroups of the right regular representation $(\Z_n)_R$ in $\sym(\Z_n)$,
and  Schur rings (S-rings for short) over $\Z_n$. It turns out that each circulant graph $\Gamma$ 
generates a suitable S-ring $\A$, such that $\aut(\Gamma)$ coincides with $\aut(\A)$.  If in addition
$\Gamma$ is a rational circulant graph, then the corresponding S-ring $\A$ is also rational.

Rational S-rings over cyclic groups were classified by Muzychuk in \cite{Muz93}.
Therefore, in principle, knowledge of \cite{Muz93} is enough in order to deduce our main results.
Nevertheless, it is helpful and natural to interpret groups of rational circulant graphs as the automorphism
groups of orthogonal group block structures on $\Z_n$. This implies interest to results of Bailey et al. about such
groups (see \cite{Bai81,BaiPRS83,Bai96}).  Consideration of orthogonal group block structures as well as of
crested products (see \cite{BaiC05}) makes it possible to describe generalized wreath products as formulas
over the alphabet with words ``crested'', ``direct'', ``wreath'', and ``symmetric group''.
Finally, the reader will be hopefully convinced that the simultaneous use of a few relatively independent 
languages, like S-rings, lattices, association schemes, posets, orthogonal block structures in conjunction  
with suitable group theoretical concepts leads naturally to the understanding of the entire picture as well as 
to a rigorous proof of the main results. 

The rest of the paper is organized as follows. Section~3 serves as a brief introduction to S-rings, while in
sections~4 and~5 we pay attention to the particular case of rational S-rings over $\Z_n$.  We conclude these
sections by crucial Corollary 5.4, which reduces the problem to the consideration of the automorphism groups of
rational S-rings over $\Z_n$.
In Section~6 an equivalent language of block structures on $\Z_n$ is introduced.
Section~7 provides the reader an opportunity to comprehend all main ideas on a level of simple examples.
In Section~8 crested products are introduced and it is shown that their use is, in principle, enough for the
recursive description of all required groups.
In Section~9 poset block structures are linked with generalized wreath products,
while Section~10 provides a relatively self-contained detailed proof of the main
Theorem~\ref{thm:main}.

A number of interesting by-product results, which follow almost immediately from the consideration are
presented in Section~11. Finally, in Section~12 we enter to a discussion of diverse historical
links between all introduced languages and techniques, though not aiming to give a comprehensive picture of all
details.

\section{Preliminaries}

In this section we collect all basic definitions and facts needed in this paper.

\subsection{Permutation groups.} 

The group of all permutations of a set $X$ is denoted by $\sym(X)$. We let $g \in \sym(X)$ act
on the right, i.e., $x^g$ is written for the image of $x$ under action of $g$, and further we have
$x^{g_1g_2}=(x^{g_1})^{g_2}$. Let $K_R$ denote the {\em right regular representation} of $K$
acting on itself, i.e., $x^k=xk$ for all $x,k \in K$. Two permutation groups $K_1 \le \sym(X_1)$ and
$K_2 \le \sym(X_2)$ are {\em permutation isomorphic} if there is a bijection
$f \colon X_1 \to X_2$, and an isomorphism $\varphi : K_1 \to K_2$ such that,
$f(x_1^{k_1})=f(x_1)^{\varphi(k_1)}$ for all $x_1 \in X_1$, $k_1 \in K_1$.

Two operations of permutations groups will play a basic role in the sequel. The
{\em permutation direct product} $K_1 \times K_2$ of groups $K_i \le \sym(X_i)$, $i=1,2$,
is the permutation representation of $K_1 \times K_2$ on $X_1 \times X_2$ acting as:
$$(x_1,x_2)^{(k_1,k_2)}=(x_1^{k_1},x_2^{k_2}), \;
(x_1,x_2) \in X_1 \times X_2, \; (k_1,k_2) \in K_1 \times K_2.$$
Note that the direct product is commutative and associative.

Let $A \le \sym(X_1)$ and $C \le \sym(X_2)$ be two permutation groups. 
The {\em wreath product} $A \wr C$ is the subgroup of $\sym(X_1 \times X_2)$
generated by the following two groups: the {\em top group} $T$, which is a faithful permutation
representation of $A$ on $X_1 \times X_2$, acting as:
$$ 
(x_1,x_2)^a=(x_1^a,x_2) \text{ for } (x_1,x_2) \in X_1 \times X_2, \; a \in A,
$$
and the {\em base group} $B$, which is the representation of the group
$C^{X_1}$ on $X_1 \times X_2$, acting as:
$$ 
(x_1,x_2)^f = (x_1,x_2^{f(x_1)}), \; (x_1,x_2) \in X_1 \times X_2, \;
f \in C^{X_1},
$$
where $f(x_1)$ is the component (belonging to $C$) of $f$, corresponding to $x_1 \in X_1$. 
(Here $C^{X_1}$ denotes the group of all functions from $X_1$ to $C$, with group operation 
$(fg)(x_1)=f(x_1) \cdot g(x_1)$ for $x_1 \in X,$ $f,g \in C^{X_1}$.) 
The group $T$ normalizes $B$, $|B \cap T|=1$, therefore 
$\sg{B,T} = B \rtimes T$. Clearly, the group $A \wr C$ has order $|A \wr C|=|T| \cdot |B| = 
|A| \cdot |C|^{|X_1|}$.
Each element $w \in A \wr C$ admits a unique decomposition $w = t b$, where $t \in T$ and
$b \in B$. Also element $w$ may be denoted as $w = [a,f(x_1)],$
called the {\em table form} of $w$ (note that here $x_1$ is a symbol for
a variable). By definition, $(x_1,x_2)^w=(x_1^a,x_2^{f(x_1)})$.
Note that, sometimes in wreath product $A \wr C$ the group $A$ is called {\em active}, while 
$C$ {\em passive} groups. 
The wreath product is associative, but not commutative.
We remark that our notation for wreath product follows, e.g., \cite{FarKM94}, and it
has opposite direction in comparison with traditions accepted in modern group theory.

A permutation group $K \le \sym(X)$ acts canonically on $X \times X$ by letting
$(x_1,x_2)^k=(x_1^k,x_2^k)$. The corresponding orbits are called the {\em 2-orbits} of $K$,
the set of which we denote by $\orb(K)$. The {\em $2$-closure} $K^{(2)}$ of $K$ is the unique
maximal subgroup of $\sym(X)$ that has the same 2-orbits as $K$. Clearly, $K \le K^{(2)}$,
and we say that $K$ is {\em $2$-closed} if $K^{(2)}=K$.

\subsection{Cayley graphs and circulant graphs.} 

By a \emph{(directed) graph} we mean a pair $\Gamma=(X,R)$, where $X$ is a nonempty set,
and $R$ is a binary relation on $X$. In the particular case when $(x,y) \in R$ if and only if $(y,x) \in R$
for all $(x,y) \in X \times X$, $\Gamma$ is also called an
{\em undirected graph}, and then
$\{x,y\}$ is said to be an (undirected) edge of $\Gamma$, which substitutes
$\{(x,y),(y,x)\}$.
The \emph{automorphism group} $\aut(\Gamma)=\aut((X,R))$ is the group of
all permutations $g$ in $\sym(X)$ that preserve $R$,
i.e., $(x^g,y^g) \in R$ if and only if $(x,y) \in R$ for all $x,y \in X$.

The {\em adjacency matrix} $A(R)$ of the relation $R$ is the $X$-by-$X$ complex
matrix defined by $$ A(R)_{x,y}=\left\{ \begin{array}{rr} 1 &\hbox{if }(x,y) \in R \\ 0
&\hbox{otherwise}.\end{array} \right.$$
The {\em eigenvalues} of $\Gamma$ are defined to be the eigenvalues of $A(R)$, and
$\Gamma$ is called  \emph{rational} if all its eigenvalues are rational.
Note that, since the characteristic polynomial of $A(R)$ has integer coefficients and
leading coefficient $\pm 1$, if its eigenvalues are rational
numbers, then these are in fact integers.

\medskip
For a subset $S \subseteq K$, the
{\em Cayley graph} $\cay(K,S)$ over $K$ with {\em connection set} $S$ is the graph
$(X,R)$ defined by
$$ X=K, \hbox{ and } R=\big\{(x,sx) \, \mid \, x \in K, s\in S\}. $$
Two immediate observations: the graph $\cay(K,S)$ is undirected if and only if
$S=S^{-1}=\{s^{-1} \mid s \in S\}$; and the right regular
representation $K_R$ is a group of automorphisms of $\cay(K,S)$.
Cayley graphs over cyclic groups are briefly called {\em circulant graphs}.

\subsection{Schur rings.} 

Let $H$ be a group written with multiplicative notation and with identity $e$.
Denote $\Q H$ the group algebra of $H$ over the field $\Q$ of rational numbers.
The group algebra $\Q H$ consists of the formal sums
$\sum_{x\in H}a_x x$, $a_x \in \Q$, equipped with entry-wise addition
$\sum_{x\in H}a_x x+\sum_{x\in H}b_x x = \sum_{x\in H}(a_x+b_x) x,$ and multiplication
$$\sum_{x\in H}a_x x \cdot \sum_{x\in H}b_x x = \sum_{x,y\in H}(a_yb_{y^{-1}x}) x.$$
Given $\Q H$-elements $\eta_1,\ldots,\eta_r$, the subspace generated by them is
denoted by $\langle \eta_1,$ $\ldots,\eta_r\rangle$. For a subset $S \subseteq H$ the {\em simple quantity} $\un{S}$
is the $\Q H$-element $\sum_{x \in H}a_x x$ with $a_x=1$ if $x \in S$,
and $a_x=0$ otherwise (see \cite{Wie64}). We shall also write $\un{s_1,\ldots,s_k}$ for the
simple quantity $\un{\{s_1,\ldots,s_k\}}$. The {\em transposed} of $\eta=\sum_{x \in H}a_x x$
is defined as $\eta^\top=\sum_{x \in H}a_x x^{-1}.$

\medskip

A subalgebra $\A$ of $\Q H$ is called a {\em Schur ring} (for short {\em S-ring})
of {\em rank $r$} over $H$ if the following axioms hold:
\begin{enumerate}[({SR}1)]
\item $\A$ has a basis of simple quantities:
$\A=\sg{ \un{T_1},\ldots,\un{T_r}}$, $T_i \subseteq H$ for all $i \in \{1,\ldots,r\}$.
\item $T_1=\{e\}$, and $\sum_{i=1}^{r}\un{T_i}=\un{H}$.
\item For every $i \in \{1,\ldots,r\}$ there exists $j \in \{1,\ldots,r\}$ such that
$\un{T_i}^\top=\un{T_j}$.
\end{enumerate}
The simple quantities $\un{T_1},\ldots,\un{T_r}$ are called the
{\em basic quantities} of $\A$, the corresponding sets $T_1,\ldots,T_r$ the
{\em basic sets} of $\A$. We set the notation $\basic(\A)=\{T_1,\ldots,T_r\}$.

\subsection{Posets and partitions.} 

A {\em partially ordered set} (for short a {\em poset}) is a pair $(X,\preceq)$, where $X$ is a nonempty set,
and $\preceq$ is a relation on $X$ which is
reflexive, antisymmetric and transitive. We write $x \prec y$ if $x \preceq y$ but
$x \ne y$.  For a subset $L \subseteq X$ we say an element $m \in L$ is {\em maximal} in $L$
if $m \preceq l$ implies $l=m$ for all $l \in L$. Similarly, $m \in L$ is {\em minimal} in $L$ if
$l \preceq m$ implies $l=m$ for all $l \in A$.
Further, we say that $i \in X$ is the {\em infimum} of
$L$ if $i \preceq l$ for all $l \in L$, and if for some $i' \in X$ we have
$i' \preceq l$ for all $l \in L$, then $i' \preceq i$.
Similarly, we say that $s \in X$ is the {\em supremum} of
$L$ if $l \preceq s$ for all $l \in L$, and if for some $s' \in X$ we have
$l \preceq s'$ for all $l \in L$, then $s \preceq s'$.
We set the notations: $i=\bigwedge L$ and $s=\bigvee L$. The infimum (supremum, respectively)
does not always exist, but if this is the case, it is determined uniquely.

The poset $(X,\preceq)$ is called a {\em lattice} if each pair of elements in
$X$ has infimum and supremum. Then we have binary operations $x \wedge y=\wedge \{x,y\}$ and
$x \vee y=\vee \{x,y\}$. The lattice $(X,\preceq)$ is {\em distributive} if for all
$x,y,z$ in $X$,
\begin{eqnarray*}
x \wedge (y \vee z) &=& (x \wedge y) \vee (x \wedge z), \\
x \vee (y \wedge z) &=& (x \vee y) \wedge (x \vee z).
\end{eqnarray*}
If $(X,\preceq)$ is a lattice, and a subset $X' \subset X$ is closed under both
$\wedge$ and $\vee$, then $(X',\preceq)$ is also a lattice, it is called
a {\em sublattice} of $(X,\preceq)$.

\medskip

Let $F$ be a partition of a set $X$.
We denote by $R_F$ the equivalence relation corresponding to $F$, and by $A(F)$ the
the adjacency matrix $A(R_F)$. We say that two partitions $E$ and $F$ of $X$ are
{\em orthogonal} if for their adjacency matrices $A(E)A(F)=A(F)A(E)$ (see \cite[Section~6.2]{Bai04} for a 
nice discussion of this concept).
The set of all partitions of $X$ is partially ordered by the relation $\sqsubseteq$,
where $E \sqsubseteq F$ ($E$ is a {\em refinement} of $F$) if any class of $E$ is contained in a
class of $F$. The resulting poset is a lattice, where $E \wedge F$ is
the partition whose classes are the intersection of $E$-classes with $F$-classes;
and $E \vee F$ is the partition whose classes are the minimal subsets being union of $E$-classes
and $F$-classes. The smallest element in this lattice is the {\em equality partition} $E_X$,
the classes of which are the singletons; the largest is the {\em universal partition} $U_X$
consisting of only the whole set $X$.

\section{More about S-rings}

Let $H$ be a finite group written with multiplicative notation and with identity $e$.
The {\em Schur-Hadamard product $\circ$} of the group algebra $\Q H$ is defined by
$$\sum_{x\in H}a_x \, x \circ \sum_{x\in H}b_x \, x:=\sum_{x\in H}a_xb_x \, x.$$
The following alternative characterization of S-rings over $H$ is a folklore
(cf. \cite[Theorem 3.1]{MuzKP01}): a subalgebra $\A$ of $\Q H$ is an S-ring if and only if
$\un{e},\un{H} \in \A$, and $\A$ is closed with respect to $\circ$ and $^{\top}$.
By this it is easy to see that the intersection of two $S$-rings is also an S-ring,
in particular, given $\Q H$-elements $\eta_1,\ldots,\eta_k$, denote
$\sgg{\eta_1,\ldots,\eta_k}$ the S-ring defined
as the intersection of all S-rings $\A$ that $\eta_i \in \A$ for all $i \in \{1,\ldots,k\}$. For
$S \subseteq H$ we shall also write $\sgg{S}$ instead of $\sgg{\un{S}},$ calling 
$\sgg{S}$ the S-ring {\em generated} by $S$.
For two S-rings $\A$ and $\B$ over $H$, we say that $\B$ is an {\em S-subring} of
$\A$ if $\B \subseteq \A$. It can be seen that this happens exactly when every basic
set of $\B$ is written as the union of some basic sets of $\A$.

Let $\A$ be an S-ring over $H$. A subset $S \subseteq H$ (subgroup $K \le H$, respectively) is an
{\em $\A$-set} ({\em $\A$-subgroup}, respectively) if $\un{S} \in \A$ ($\un{K} \in \A$, respectively).
If $S \subseteq H$ is an $\A$-set, then $\sg{S}$ is an $\A$-subgroup (see \cite[Proposition~23.6]{Wie64}).
By definition, the trivial subgroups
$\{e\}$ and $H$ are $\A$-subgroups, and for two $\A$-subgroups $E$ and $F$, also $E \cap F$
and $\sg{E,F}$ are $\A$-subgroups. In other words, the $\A$-subgroups form a sublattice of
the subgroup lattice of $H$. Let $K$ be an $\A$-subgroup.
Define $\A_{K}=\A \cap \Q K$. It is easy to check that
$\A_K$ is an S-ring over $K$ and
$$ \basic(\A_K)=\{ T \in \basic(\A) \, \mid \, T \subseteq K\}.$$
We shall call $\A_K$ an {\em induced S-subring} of $\A$.

Following \cite{KliP81}, by an {\em automorphism} of an S-ring
$\A = \sg{\un{T}_1,\ldots,\un{T}_r}$ over $H$ we mean a permutation $f \in \sym(H)$ which is an
automorphism of all {\em basic graphs} $\cay(H,T_i)$. Thus the automorphism group of $\A$ is
$$ \aut(\A) = \bigcap_{i=1}^{r} \aut(\cay(H,T_i)).$$

\medskip

The simplest examples of an S-ring are the whole group algebra $\Q H$, and
the subspace $\sg{\un{e},\un{H \setminus \{e\}}}$. The latter is called the
{\em trivial S-ring} over $H$. Further examples are provided by permutation
groups $G$ which are overgroups of $H_R$ in $\sym(H)$ (i.e., $H_R \le G \le \sym(H)$).
Namely, letting $T_1=\{e\}, T_2,\ldots,T_r$ be the orbits of the stabilizer $G_e$ of $e$ in $G$,
it follows that the subspace $\sg{\un{T}_1,\ldots,\un{T}_r}$ is an S-ring over $H$
(see \cite[Theorem 24.1]{Wie64}). This fact was proved by Schur, and the resulting S-ring
is also called the {\em transitivity module} over $H$ induced
by the group $G_e$, notation $V(H,G_e)$. It turns out that not every S-ring over
$H$ arises in this way, and we call therefore an S-ring $\A$ {\em Schurian} if
$\A=V(H,G_e)$ for a suitable overgroup $G$ of $H_R$ in $\sym(H)$.
The connection between permutation groups and S-rings is reflected in the
following proposition (see \cite[Theorem 3.13]{MuzKP01}).

\begin{prop}
Let $\A$ and $\B$ be arbitrary S-rings over $H$, and let
$G$ and $K$ be arbitrary overgroups of $H_R$ in $\sym(H)$. Then
\begin{enumerate}[(i)]
\item $\A \subseteq \B \Rightarrow \aut(\A) \ge \aut(\B)$.
\item $G \le K \Rightarrow V(H,G_e) \supseteq V(H,K_e)$.
\item $\A \subseteq V(H,\aut(\A)_e)$.
\item $G \le \aut(V(H,G_e))$.
\end{enumerate}
\end{prop}

\noindent
The above proposition describes a {\em Galois correspondence} between S-rings over $H$ and
overgroups of $H_R$ in $\sym(H)$. We remark that it is a particular case of a
Galois correspondence between coherent configurations and permutation groups (cf. \cite{Wei76,FarIK90}).

\medskip

The starting point of our approach toward Theorem~\ref{thm:main} is the following consequence of the
Galois correspondence, which is formulated implicitly in \cite{Wei76}.

\begin{thm}\label{thm:W}
Let $H$ be a finite group and $S \subseteq H$. Then
$$ \aut(\cay(H,S)) = \aut(\sgg{S}). $$
\end{thm}

\section{Rational S-rings over cyclic groups}

In this section we turn to S-rings over cyclic groups.
Our goal is to provide a description of those S-rings $\A$ that
$\A=\sgg{S}$ for some rational circulant graph $\cay(\Z_n,S)$.

Throughout the paper the cyclic group of order $n$ is given by the additive
cyclic group $\Z_n$, written as $\Z_n=\{0,1,\ldots,n-1\}$.
Note that, we have switched from multiplicative to additive notation.
For a positive divisor $d$ of $n$, $Z_d$ denotes the unique subgroup of $\Z_n$
of order $d$, i.e.,
$$ Z_d =\sg{m}=\big\{ xm \, \mid \, x \in \{0,\ldots,d-1\} \big\}, \textrm{ where } n=dm.$$
Let $\Z_n^*=\{ i \in \Z_n \, \mid \, \gcd(i,n)=1\}$, i.e., the multiplicative group
of invertible elements in the ring $\Z_n$.
(By some abuse of notation $\Z_n$ stands parallel for both the ring and also its
additive group. Moreover, sometimes $\Z_n$ simply denotes the set of symbols $0,1,\ldots,n-1$.) 
For $m \in \Z_n^*$, and a subset
$S \subseteq \Z_n$, define $S^{(m)}=\{ ms  \mid s \in S \}$. Two subsets
$R,S \subseteq \Z_n$ are said to be {\em conjugate} if
$S=R^{(m)}$ for some $m \in \Z_n^*$. The {\em trace} $\trS$ of $S$ is the union of all subsets conjugate to $S$,
i.e., 
$$ \trS = \bigcup_{m\in \Z_n^*}S^{(m)}.$$
The elements $m$ in $\Z_n^*$ act on $\Z_n$ as group of automorphisms by sending $x$ to $mx$.
We have corresponding orbits
\begin{equation}\label{eq:Znd}
(\Z_n)_d = \big\{ x \in \Z_n \, \mid \,  \gcd(x,n)=d \big\},
\end{equation}
where $d$ runs over the set of positive divisors of $n$.
The {\em complete S-ring of traces} is the transitivity module
$$ V(\Z_n,\Z_n^*)=\sg{ \un{(\Z_n)_d} \, \mid \, d \mid n }.$$
By the {\em rational} (or {\em trace})
S-rings over $\Z_n$ we mean the S-subrings of $V(\Z_n,\Z_n^*)$.
For an S-ring $\A$ over $\Z_n$ its {\em rational closure} $\trA$ is the
S-ring defined as $\trA=\A \cap V(\Z_n,\Z_n^*)$, and thus $\A$ is rational if and only if $\A=\trA$.

\medskip

Recall that a circulant graph $\cay(\Z_n,S)$ is rational if it has a rational spectrum.
The following result describes its connection set $S$ in terms of the generated S-rings $\sgg{S}$
(cf. \cite{BriM79}).

\begin{thm}\label{thm:BM}
A circulant graph $\Gamma=\cay(\Z_n,S)$ is rational if and only if the generated S-ring
$\sgg{S}$ is a rational S-ring over $\Z_n$.
\end{thm}

It follows from the theorem that $S$ is a union of some sets of the form $(\Z_n)_d$.
In particular, exactly $2^{\tau(n)-1}$ subsets of $\Z_n$ define a rational circulant
graph without loops (i.e., $0 \notin S$). Here $\tau(n)$ denotes the number of positive 
divisors of $n$. As we shall see in 11.1, the resulting graphs are pairwise
non-isomorphic.

\section{Properties of rational S-rings over cyclic groups}

Denote $L(n)$ the lattice of positive divisors
of $n$ endowed with the relation $x \mid y$ ($x$ divides $y$).
For two divisors $x$ and $y$, we write $x \wedge y$ for their greatest common
divisor, and $x \vee y$ for their least common multiple. Note that, the lattice $L(n)$ is distributive,
and if $L$ is any set of positive divisors of $n$, then the poset $(L,\mid)$
is a sublattice of $L(n)$ if and only if $L$ is closed with respect to $\wedge$ and $\vee$.
By some abuse of notation this sublattice we shall denote by $L$ as well.

For a sublattice $L$ of $L(n)$, and $m \in L$, we define the sets
$$ 
L\subl{m} = \{ x \in L \, \mid \, x \mid m \}, \text{ and }
L\supl{m} = \{ x \in L \, \mid \, m \mid x \}.
$$
It is not hard to see that these are sublattices of $L(n)$.

\medskip

The following classification of rational S-rings over $\Z_n$ is due to Muzychuk
(see \cite[Main Theorem]{Muz93}).

\begin{thm}\label{thm:M} $\,$
\begin{enumerate}[(i)]
\item Let $L$ be a sublattice of $L(n)$ such that $1,n \in L$. Then
the vector space $\A=\sg{\un{Z_l} \, \mid \, l \in L}$ is an S-ring over $\Z_n,$ which is rational.
\item Let $\A$ be a rational S-ring over $\Z_n$. Then there exists a sublattice
$L$ of $L(n)$, $1,n \in L$, such that $\A=\sg{\un{Z_l} \, \mid \, l \in L}$.
\end{enumerate}
\end{thm}

We remark that if $\A=\sg{\un{Z_l} \, \mid \, l \in L}$ is the S-ring in part (i) above, then  
the simple quantities $\un{Z_l}$ form a basis of the vector space $\A$, where $l$ runs over the set $L$. 
This basis we shall also call the  {\em group basis} of $\A$. It is also true that all $\A$-subgroups 
appear in this basis, i.e., for any subgroup $Z_k \le \Z_n,$ we have $\un{Z_k} \in \A$ if and only if $k \in L$. 
The basic quantities of the rational S-ring $\A$ are easily obtained from its group basis, 
namely $\basic(\A)$ consists of the sets:
\begin{equation}\label{eq:whatZl}
\widehat{Z}_l = Z_l \setminus \bigcup_{d \in L\subl{l},d<l}Z_d, \quad l\in L.
\end{equation}

\medskip

In the rest of this section we are going to prove that rational S-rings over 
$\Z_n$ are generated by subsets of $\Z_n$. More formally, that every rational S-ring $\A$ over $\Z_n$ satisfies 
$\A=\sgg{S},$ where $S$ is a suitable subset $S \subseteq \A$. Notice that, the corresponding circulant 
graph $\cay(\Z_n,S)$ is rational (see Theorem~\ref{thm:BM}), and 
its automorphism group $\aut(\cay(\Z_n,S)) = \aut(\A)$ (see Theorem~\ref{thm:W}). 

\medskip 

We start with an auxiliary lemma, for which the authors thank 
Muzychuk (see \cite{Muz}).

\begin{lem}\label{lem:M}
Let $L$ be a sublattice of $L(n)$, $1,n\in L$. Let $m$ be a maximal element of
the poset $(L \setminus \{n\},\mid)$, and $s$ be the
smallest number in the set $L \setminus L\subl{m}$. Then
$$ 
L \setminus L\subl{m} = \Big\{  x \, \frac{s}{m \wedge s} \, \mid \,
x \in (L\subl{m})\supl{m \wedge s}  \Big\}.
$$
\end{lem}

\proof Define the mapping
$$ 
f \colon L \setminus L\subl{m} \to L\subl{m}, \; l \mapsto m \wedge l. 
$$
Let $l \in L \setminus L\subl{m}$. As $m$ is maximal, $l \vee m=s \vee m=n$. By distributive law,
$(l \wedge s) \vee m=(l \vee m) \wedge (s \vee m)=n$.
Thus $l \wedge s \in L \setminus L\subl{m}$, and by the choice of $s$, $s \le l \wedge s$,
hence $s \mid l$, $(m \wedge s) \mid f(l)$, and $f(l) \in (L\subl{m})\supl{m \wedge s}$.

On the other hand, choose $x \in (L\subl{m})\supl{m \wedge s}$.
Then $l=s \vee x$ is in $L \setminus L\subl{m}$, and we find $f(l)=m \wedge l=
(m \wedge s) \vee (m \wedge x)=(m \wedge s) \vee x=x$. Also,
$f(L \setminus L\subl{m}) = (L\subl{m})\supl{m \wedge s}$.

For each $l \in L \setminus L\subl{m}$,
\begin{equation}\label{eq in lem:M}
s \vee f(l) = s \vee (m \wedge l)=(s \vee m) \wedge (s \vee l)=n \wedge l=l.
\end{equation}
The lemma follows as
$$ 
L \setminus L\subl{m} = \Big\{ s \vee f(l) \, \mid \, l \in L  \Big\} =
\Big\{  s \vee x =  \frac{s}{m \wedge s} \, x \, \mid \, x \in (L\subl{m})\supl{m \wedge s} \Big\},
$$
here we use the property $x \wedge s=m \wedge s$. \QED

\begin{prop}\label{prop:convBM}
Let $\A$ be a rational S-ring over $\Z_n$. Then there exists a
subset $S \subseteq \Z_n$ such that $\A=\sgg{\un{S}}$.
\end{prop}

\proof We proceed by induction on $n$. The case $n=1$ is trivially true. 
Let $n > 1$. By (ii) of Theorem~\ref{thm:M},
\begin{equation}\label{eq1 in prop:convBM}
\A=\sg{ \un{Z_l} \, \mid \, l \in L},
\end{equation}
where $L$ is a sublattice of $L(n)$, $1,n \in L$. Let $m$ be a maximal element in the poset
$(L \setminus \{n\},\mid)$, and $s$ be the smallest number in the set $L \setminus L\subl{m}$.
Apply the induction hypothesis to the induced S-subring $\A|_{Z_m}=\A \cap \Q Z_m$. This results in a
subset $R \subseteq Z_m$ such that $\A|_{Z_m}=\sgg{R}$.
Pick the basic set $\widehat{Z}_s \in \basic(\A)$, see \eqref{eq:whatZl}. 
By the choice of $s$ we get 
$$ 
\widehat{Z}_s = Z_s \setminus \bigcup_{d \in L\subl{s}, d < s}Z_d = 
Z_s \setminus Z_{m \wedge s}. 
$$ 
Let 
$$ 
S = R \cup \widehat{Z}_s, \text{ and } \A'=\sgg{S}.
$$
It is clear that $S$ equals its trace $\trS$, so $\A'$ is a rational S-ring.
We complete the proof by showing that in fact $\A=\A'$.

As $\un{\widehat{Z}_s} \in \A$, $\un{S} \in \A$,
hence $\A' \subseteq \A$. By \eqref{eq1 in prop:convBM}, to have $\A \subseteq \A'$
it is enough to show that, for any positive divisor $l$ of $n$,
\begin{equation}\label{eq2 in prop:convBM}
l \in L \implies \un{Z_l} \in \A'. 
\end{equation}

We show first that $\un{Z_s} \in \A'$. 
Let $T \in \basic(\A')$ such that $(\Z_n)_{n/s} \subseteq T$.
Consider the subgroup $\sg{T}$, and let $\sg{T}=Z_t$. 
As $\A' \subseteq \A,$ $\un{T} \in \A,$ and therefore $\un{Z_t}$ is in $\A$. This gives $t \in L$. 
Clearly, $t \in L \setminus L\subl{m},$ and hence $t = s \vee (m \cap t),$ see \eqref{eq in lem:M}. 
It follows from the description of basic sets in \eqref{eq:whatZl} that 
$T$ contains a generator of $\sg{T}=Z_t$. Thus if $t \ne s,$ then 
$T \cap (\Z_n \setminus Z_m \setminus Z_s) \ne \emptyset$. But, 
$T \subseteq S$ and $S \subseteq Z_m \cup Z_s,$ implying that $t=s,$ and so $\un{Z_s}$ is in $\A'$.  

Thus $\un{S \setminus Z_s}=\un{R \setminus Z_s} \in \A'$. 
Let $s < n$. We may further assumed that $R \cap (\Z_n)_{n/m} \ne \emptyset,$ 
otherwise replace $R$ with its complement in $Z_m \setminus \{0\}$. 
Thus we find $\un{Z_m}=\un{\sg{R \setminus Z_s}} \in \A'$. If $s = n,$ then 
we may assume that $(Z_m \setminus R) \cap (\Z_n)_{n/m} \ne \emptyset$. From this  
$\un{Z_m}=\un{\sg{\Z_n \setminus S}} \in \A'$. Then 
$$ \A|_{Z_m}=\sgg{\un{R}} \subseteq \A'|_{Z_m} \subseteq \A|_{Z_m},$$
from which $\A|_{Z_m}=\A'|_{Z_m}$. We conclude that \eqref{eq2 in prop:convBM} holds if $l \in L\subl{m}$.

Let $l \in L \setminus L\subl{m}$. By \eqref{eq in lem:M} we can
write $l = s \vee l',$ where $l' =  m \wedge l$ is in $L\subl{m}$. 
Then $Z_l=\sg{ Z_{l'}, Z_s}$. As both $\un{Z_{l'}} \in \A'$ and $\un{Z_s} \in \A'$, $\un{Z_l} \in \A'$ follows, 
and this completes the proof of \eqref{eq2 in prop:convBM}. \QED

By Theorems~\ref{thm:W},~\ref{thm:BM} and Proposition~\ref{prop:convBM},
we obtain the following equivalence.

\begin{cor}\label{cor:step1}
Let $G$ be a permutation group acting on the cyclic group $\Z_n$.
The following are equivalent:
\begin{enumerate}[(i)]
\item $G=\aut(\cay(\Z_n,S))$ for some rational circulant graph $\cay(\Z_n,S)$.
\item $G=\aut(\A)$ for some rational S-ring $\A$ over $\Z_n$.
\end{enumerate}
\end{cor}

\section{From rational S-rings to block (partition) structures}

A {\em block structure} $\F$ on a set $X$ is simply a collection of partitions of $X$.
A partition $F$ of $X$ is {\em uniform} if all classes of $F$ are of the same
cardinality. Block structure $\F$ is called {\em orthogonal} (see e.g. \cite{Bai04}) if 
the following axioms hold:

\begin{enumerate}[({OBS}1)]
\item $E_X,U_X \in \F$.
\item Every $F \in \F$ is uniform.
\item Every two $E,F \in \F$ are orthogonal.
\item For every two $E,F \in \F$, both $E \wedge F \in \F$ and $E \vee F \in \F$.
\end{enumerate}
Note that, if $\F$ is orthogonal, then the poset $(\F,\sqsubseteq)$ is a
lattice, where $\sqsubseteq$ is the refinement relation defined on the set of partitions of 
$X$. Below we say that $\F$ is {\em distributive} if the lattice
$(\F,\sqsubseteq)$ is distributive.

\medskip

The following example of a block structure is crucial in the sequel.

\begin{exm} ({\em group block structure}) \
Let $H$ be an arbitrary group, and $K$ be a subgroup of $H$. Denote by $F_K$ the partition
of $H$ into right cosets of $K$. A {\em group block structure} on $H$ is a block structure
$(H,\{F_K \, \mid \, K \in \mathcal{K}\})$ where $\mathcal{K}$ is a set of subgroups of
$H$ satisfying the following axioms:

\begin{enumerate}[({GBS}1)]
\item The trivial subgroup $\{e\}$ is in $\mathcal{K}$.
\item For every two $K_1,K_2 \in \mathcal{K}$, $K_1K_2=K_2K_1$, and $K_1K_2 \in \mathcal{K}$.
\end{enumerate}

\noindent 
It follows that the group block structure $(H,\{F_K \, \mid \, K \in \mathcal{K}\})$ is
orthogonal if and only if $H \in \mathcal{K}$, and $(\mathcal{K},\le)$ is a sublattice of the
subgroup lattice of $H$. \QED
\end{exm}

In this context Theorem~\ref{thm:M} can be rephrased as follows.

\begin{thm}\label{thm:M2} $\,$
\begin{enumerate}[(i)]
\item Let $\F$ be an orthogonal group block structure on $\Z_n$. Then the
vector space $\A=\sg{ \un{Z_l} \, \mid \, F_{Z_l} \in \F }$ is an S-ring over $\Z_n$.
\item Let $\A$ be a rational S-ring over $\Z_n$. Then there exists
an orthogonal group block structure $\F$ on $\Z_n$ such that
$\A=\sg{ \un{Z_l} \, \mid \, F_{Z_l} \in \F }$ (here again equality means equality of vector spaces).
\end{enumerate}
\end{thm}

\medskip

For $i=1,2$, let $\F_i$ be a block structure on $X_i$. 
Following \cite{Bai96}, a {\em weak isomorphism} from $\F_1$ to $\F_2$ is a 
bijection $f \colon X_1 \to X_2$ such that there exists a bijection
$g \colon \F_1 \to \F_2$ for which $(x_1,y_1) \in  R_F$ if and only if
$(x_1^f,y_1^f) \in R_{F^g}$ for all $x_1,y_1 \in X_1$, and $F \in \F_1$.
The mapping $f$ is also called a {\em strong isomorphism with respect to $g$}, or
simply a {\em strong isomorphism} if $g$ is understood.
In particular, a {\em weak automorphism} of
$\F$ is a weak isomorphism of $\F$ onto itself, and a {\em strong automorphism} (or an
{\em automorphism}) is a weak automorphism which is strong with respect to the identity.
The {\em automorphism group} $\aut(\F)$ of $\F$ is therefore the permutation group (see also \cite{Bai81})
$$ \aut(\F)=\bigcap_{F \in \F} \aut((X,R_F)).$$


\begin{prop}\label{prop:autA=autF}
Let $\A$ be a rational S-ring over $\Z_n$, and $\F$ be an orthogonal group block structure
on $\Z_n$ such that $\A=\sg{ \un{Z_l} \, \mid \, F_{Z_l} \in \F }$. Then
$\aut(\A)=\aut(\F)$.
\end{prop}

\proof Let $L$ be the sublattice of $L(n)$ corresponding to $\F$.
To ease notation, we write $R_l$ for the
relation $R_{F_{Z_l}},$ where $l \in L$. Then $\A$ has basic sets $\widehat{Z}_l$, $l\in L$,
see \eqref{eq:whatZl}. Let $\widehat{R}_l$ be the relation on $\Z_n$ that is given by 
the arc set of $\cay(\Z_n,\widehat{Z}_l)$, i. e., 
$\cay(\Z_n,\widehat{Z}_l)=(\Z_n,\widehat{R}_l)$. Thus for $l \in L$,
$$ 
\widehat{R}_l=R_l \setminus \bigcup_{d \in L\subl{l},d<l}R_d, \text{ and }
R_l=\bigcup_{d \in L\subl{l}}\widehat{R}_d.
$$
Thus for $g \in \aut(\A)$, $R_l^g=\cup_{d \in L\subl{l}}\widehat{R}_d^{\; g}=
\cup_{d \in L\subl{l}}\widehat{R}_d=R_l$, and so $g \in \aut(\F)$.
Similarly, if $g \in \aut(\F)$, then
$\widehat{R}_l^{\, g}=R_l^g \setminus \cup_{d \in L\subl{l},d< l}R_d^g=
R_l \setminus \cup_{d \in L\subl{l},d< l}R_d=\widehat{R}_l$, implying
$g \in \aut(\A)$. Therefore $\aut(\A)=\aut(\F)$. \QED

\noindent
We remark that the above correspondence in Theorem~\ref{thm:M2} is a particular case of a
correspondence between orthogonal block structures and association
schemes, see the discussion in 11.2.

\medskip

By Corollary~\ref{cor:step1}, Theorem~\ref{thm:M2}, and Proposition~\ref{prop:autA=autF},
we obtain the following equivalence.

\begin{cor}\label{cor:step2}
Let $G$ be a permutation group acting on the cyclic group $\Z_n$.
The following are equivalent:
\begin{enumerate}[(i)]
\item $G=\aut(\cay(\Z_n,S))$ for some rational circulant graph $\cay(\Z_n,S)$.
\item $G=\aut(\F)$ for some orthogonal group block structure $\F$ on $\Z_n$.
\end{enumerate}
\end{cor}

\section{Simple examples}

We interrupt the main line of the presentation, exposing a few simple examples.  The goal is to provide the
reader additional helpful context. 
Recall that according to the previous propositions each
rational S-ring over $\Z_n$ is uniquely determined by a suitable sublattice of the lattice $L(n)$, or in equivalent terms, by a suitable block structure on $\Z_n$.  Moreover, for each rational S-ring a Cayley graph may be found
which generates the S-ring in certain prescribed sense.  Nevertheless, in many cases consideration of several Cayley graphs in role of generators allows to better comprehend the considered S-ring.  Each time in this section we abuse notation, identifying lattices with their S-rings.

\medskip

Our first example refines Example~\ref{ex:n=6}.

\begin{exm}\label{ex:n=6'}
Here $n=6$, we first depict lattice $L=L(6)$.
Clearly $L$ has $3$ sublattices containing $1$ and $6$ as shown in Figure~2.
$\aut(L_0) = S_6$.
The sublattice $L_1$ is generated by the point $3$, which may be regarded as partition
$\{\{0,2,4\},\{1,3,5\}\}$.  $\aut(L_1)$ is recognized as the wreath product
$S_2 \wr S_3$ of order $2 \cdot (3!)^2 = 72$.
Similarly, $\aut(L_2)$ is the wreath product of order $3! \cdot (2!)^3= 48$.
A significant message is that, for lattice $L$ we have
$\aut(L) = \aut(L_1) \cap \aut(L_2)= S_3 \times S_2$, a transitive group of order $12$,
containing $(\Z_6)_R$ as a subgroup. \QED
\end{exm}

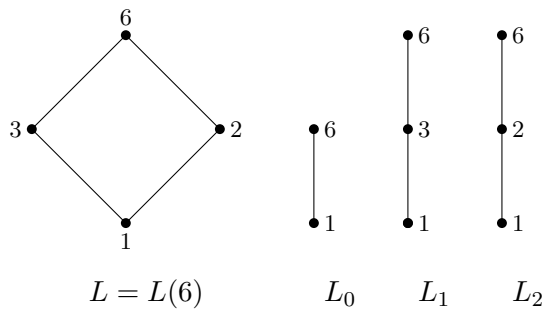
\begin{figure} [h!]
\centering
\begin{tikzpicture}[scale=1.25] ---
\foreach \i in {0}
{\draw (\i ,0) -- (\i +1,1) -- (\i+2 ,0) -- (\i +1,-1) -- cycle;
\fill (\i, 0) circle (1.5pt);
\fill (\i +1,1) circle (1.5pt);
\fill (\i+2 ,0) circle (1.5pt);
\fill (\i +1,-1) circle (1.5pt);}
\draw (0,0) node[left] {\footnotesize $3$};
\draw (1,1) node[above] {\footnotesize $6$};
\draw (2,0) node[right] {\footnotesize $2$};
\draw (1,-1) node[below] {\footnotesize $1$};
{\draw (3,-1) -- (3,0);}
\fill (3,-1) circle (1.5pt);
\fill (3,0) circle (1.5pt);
\draw (3,-1) node[right] {\footnotesize $1$};
\draw (3,0) node[right] {\footnotesize $6$};
{\draw (4,-1) -- (4,0) -- (4,1);}
\fill (4,-1) circle (1.5pt);
\draw (4,1) node[right] {\footnotesize $6$};
\draw (4,0) node[right] {\footnotesize $3$};
\draw (4,-1) node[right] {\footnotesize $1$};
{\draw (5,-1) -- (5,0) -- (5,1);}
\draw (5,1) node[right] {\footnotesize $6$};
\draw (5,0) node[right] {\footnotesize $2$};
\draw (5,-1) node[right] {\footnotesize $1$};
\draw (0.5,-1.75) node[right] {\normalsize $L=L(6)$};
\draw (3,-1.75) node[right] {\normalsize $L_0$};
\draw (4,-1.75) node[right] {\normalsize $L_1$};
\draw (5,-1.75) node[right] {\normalsize $L_2$};
\foreach \k in {-1,0,1}
{
\fill (4,\k) circle (1.5pt);
\fill (5,\k) circle (1.5pt);
}
\end{tikzpicture}
\caption{Sublattices of $L(6)$.}
\end{figure}

Next two rules appear as natural generalization of the observations learned from Example~\ref{ex:n=6'}.
Recall that a partition $E$ is a refinement of partition $F$ if 
each class of $E$ is a part of some class of $F$. 

\medskip

\noindent{\bf Rule 1.} \ The partition defined by node $k$ is a refinement of the 
partition defined by node $kl,$ see part (i) of Figure~3.
This is also called {\em nesting} of  partitions (see \cite{BaiC05}).
In this case $\aut(L) = S_l \wr S_k$.

\medskip

\noindent{\bf Rule 2.} \ 
Let $\gcd(k,l)=1$.
Each class of the partition defined by node $kl$ is union of
classes defined by nodes $k$ and $l$, respectively, such that the latter partitions have classes 
intersecting in at most one element, see part (ii) of Figure~3.
This is also called {\em crossing} of  partitions (see \cite{BaiC05}).
In this case $\aut(L) = S_k \times S_l$.

\medskip

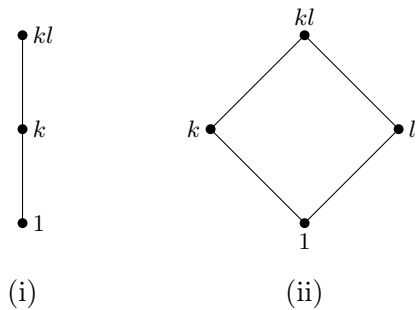
\begin{figure}[h!]
\centering
\begin{tikzpicture}[scale=1.25] ---
{\draw (0,0) -- (0,1) -- (0,2);}
\fill (0,0) circle (1.5pt);
\fill (0,1) circle (1.5pt);
\fill (0,2) circle (1.5pt);
\draw (0,0) node[right] {\footnotesize $1$};
\draw (0,1) node[right] {\footnotesize $k$};
\draw (0,2) node[right] {\footnotesize $kl$};
\draw (0,-0.75) node {(i)};

\draw (3,0) -- (4,1) -- (3,2) -- (2,1) -- (3,0);
\fill (3,0) circle (1.5pt);
\fill (4,1) circle (1.5pt);
\fill (3,2) circle (1.5pt);
\fill (2,1) circle (1.5pt);
\draw (3,0) node[below] {\footnotesize $1$};
\draw (4,1) node[right] {\footnotesize $l$};
\draw (3,2) node[above] {\footnotesize $kl$};
\draw (2,1) node[left] {\footnotesize $k$};
\draw (3,-0.75) node {(ii)};
\end{tikzpicture}
\caption{Rules $1$ and $2$.}
\end{figure}

The following simple reductions rules are 
clear generalizations of the above Rules 1 and 2.

\medskip

\noindent{\bf Reduction rule 1.} \  
This falls into two cases: either each partition defined by node $i, i \neq lm,$ is a refinement of the partition defined by node $m$, see (i) of Figure~4; or the partition defined by node $l$ is a refinement of 
each partition defined by node $i, i \neq 1,$ see (ii) of Figure~4. 
In the first case $\aut(L) = S_l \wr \aut(L_1)$, and in the second case $\aut(L) = \aut(L_1) \wr S_l$.

\medskip 

\noindent{\bf Reduction rule 2.} \ 
Here $n=ij, \gcd(i,j)=1,$ and $L = L_1 \times L_2$ is a {\em direct product} of sublattices $L_1$ of $L(i)$ and 
$L_2$ of $L(j)$. An essential property of such situation is that the entire lattice $L$ contains a 
sublattice, isomorphic to (ii) in Figure~3. (This fact is conditionally depicted in (iii) of Figure~4.)
In this case $\aut(L) = \aut(L_1) \times \aut(L_2)$.
 

\begin{figure} [h!]
\centering

\begin{tikzpicture}[scale=1.5] ---
\draw (0,0) circle (0.5cm);
\draw (0,0) node {\normalsize $L_1$};
\fill (0,0.5) circle (1pt);
\draw (0,0.6) node[right] {\footnotesize $m$};
\fill (0,1.5) circle (1pt);
\draw (0,1.5) node[right] {\footnotesize $lm$};
\draw (0,0.5) -- (0,1.5);
\draw (0,-1) node {(i)};

\draw (2,1) circle (0.5cm);
\draw (2,1) node {\normalsize $L_1$};
\fill (2,-0.5) circle (1pt);
\draw (2,-0.5) node[right] {\footnotesize $1$};
\fill (2,0.5) circle (1pt);
\draw (2,0.4) node[right] {\footnotesize $l$};
\draw (2,-0.5) -- (2,0.5);
\draw (2,-1) node {(ii)};

\draw (4,0.5) circle (0.5cm);
\draw (4,0.5) node {\normalsize $L_1$};
\draw (5.5,0.5) circle (0.5cm);
\draw (5.5,0.5) node {\normalsize $L_2$};
\fill (4.75,-0.5) circle (1pt);
\draw (4.75,-0.5) node[right] {\footnotesize $1$};
\fill (4.75,1.5) circle (1pt);
\draw (4.75,1.5) node[right] {\footnotesize $ij$};
\fill (4.23,0.92) circle (1pt);
\draw (4.23,0.92) node[above] {\footnotesize $i$};
\fill (5.27,0.92) circle (1pt);
\draw (5.27,0.92) node[above] {\footnotesize $j$};
\draw (4.75,-0.5) -- (4.23,0.08); 
\draw (4.75,-0.5) -- (5.27,0.08);
\draw (4.23,0.92) -- (4.75,1.5);
\draw (5.27,0.92) -- (4.75,1.5);
\draw (4.75,-1) node {(iii)};

\end{tikzpicture}
\caption{Reduction rules 1 and 2.} 
\end{figure}
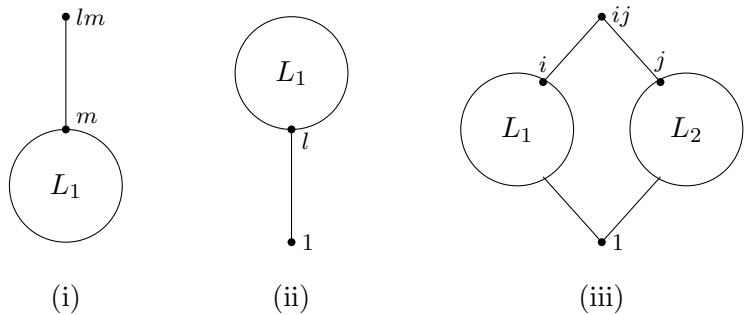

\medskip 

The created small toolkit of rules proves immediately its efficiency.

\begin{exm}\label{ex:=pe}
Here $n=p^e$, $p$ is a prime number.
In this case each sublattice of $L(p^e)$ forms a chain, hence can be constructed with 
using only Reduction rule 1. Thus the automorphism group of each sublattice of $L(p^e)$ is an 
iterated wreath product of symmetric groups. \QED
\end{exm}

\begin{exm}\label{ex:n=pqr}
Here $n= pqr$, $p,q$, and $r$ are distinct primes.
One can case by case describe possible sublattices of $L(pqr)$ and in each case to express corresponding
automorphism group with the aid of operation of direct and wreath products.

\begin{figure} [h!]
\centering

\begin{tikzpicture}[scale=1.5] ---
\draw (0,0) -- (1,1) -- (0,2) -- (-1,1) -- (0,0);
\draw[style=dotted] (1,2) -- (0,1) -- (-1,2) -- (0,3);
\draw[style=dotted] (0,0) -- (0,1);
\draw (1,2) -- (0,3);
\draw (1,1) -- (1,2);
\draw[style=dotted] (-1,1) -- (-1,2);
\draw (0,2) -- (0,3);
\fill (0,0) circle (1.5pt);
\foreach \i in {-1,0,1}
{\fill (\i,1) circle (1.5pt);
\fill (\i,2) circle (1.5pt);}
\fill (0,3) circle (1.5pt);
\draw (0,0) node[right] {\footnotesize $1$};
\draw (-1,1) node[right] {\footnotesize $p$};
\draw (0,1) node[right] {\footnotesize $q$};
\draw (1,1) node[right] {\footnotesize $r$};
\draw (-1,2) node[right] {\footnotesize $pq$};
\draw (0,2) node[right] {\footnotesize $pr$};
\draw (1,2) node[right] {\footnotesize $qr$};
\draw (0,3) node[right] {\footnotesize $pqr$};
\end{tikzpicture}

\caption{Sublattice $L$ of $L(pqr)$.}
\end{figure}
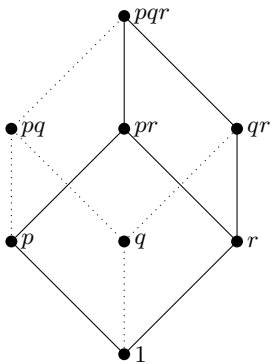

For example, for the sublattice $L$ in Figure~5 we easily obtain
$\aut(L) = (S_q \wr S_r) \times S_p.$  
(Indeed, here $L$ is a direct product of two chains with $2$ and $3$ nodes.)

\QED
\end{exm}

It is not true however that such easy life is possible for arbitrary value of $n$.
A simple case of a failure is provided by $n=p^{2}q^{2}$, where $p,q$ are distinct primes.
To make presentation more clear and visible let us consider a concrete sublattice $L$ of $L(36)$.

\begin{exm}\label{ex:n=36}
Here $n=36$, and let $L$ be the sublattice of $L(36)$ given in Figure~6.

\begin{figure}[h!]
\centering

\begin{tikzpicture}[scale=1.25] ---
\fill (0,0) circle (1.5pt);
\draw (0,0) node[right] {\footnotesize $1$};
\fill (-1,1) circle (1.5pt);
\draw (-1,1) node[right] {\footnotesize $3$};
\fill (1,1) circle (1.5pt);
\draw (1,1) node[right] {\footnotesize $2$};
\fill (0,2) circle (1.5pt);
\draw (0,2) node[right] {\footnotesize $6$};
\fill (2,2) circle (1.5pt);
\draw (2,2) node[right] {\footnotesize $4$};
\fill (1,3) circle (1.5pt);
\draw (1,3) node[right] {\footnotesize $12$};
\fill (-1,3) circle (1.5pt);
\draw (-1,3) node[right] {\footnotesize $18$};
\fill (0,4) circle (1.5pt);
\draw (0,4) node[right] {\footnotesize $36$};
\draw (0,0) -- (-1,1) -- (0,2) -- (1,1) -- (0,0);
\draw (1,1) -- (2,2) -- (1,3) -- (0,2);
\draw (1,3) -- (0,4) -- (-1,3) -- (0,2);
\end{tikzpicture}

\caption{Sublattice $L$ of $L(36)$.}
\end{figure}
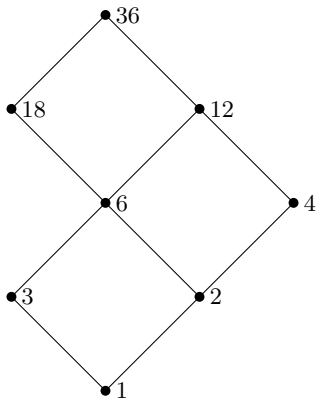

At this stage we wish to describe the automorphism group of $L$, using transparent, clear and intuitive arguments, avoiding however more rigorous justification. We note that we will return to the lattice
$L$ in this example a few times in our further presentation.
It may be convenient for us to identify the group $\aut(L)$ with the group $\aut(\Gamma)$, for a suitable
Cayley graph $\Gamma$.  Recall that as a rule, one may find several possibilities to reach such graph
(cf. Section 5). We however wish to use first a more dogmatic (in a sense naive) approach, which is based
completely on the paper \cite{Muz94}.
Basing on this text, we easily identify the unique S-ring which corresponds to $L$.
(We admit that our theoretical reasonings were, in addition, confirmed independently with the aid of a
computer via the use of COCO (see \cite{FarK91}).)
Thus we reach that the S-ring defined by $L$ has rank $8$ with the basic sets $B_k$ as follows
(see also \eqref{eq:whatZl}):
$$ Q_{36}, \; Q_2 \cup Q_4, \; Q_3, \; Q_6, \; Q_9, \; Q_{12}, \; Q_{18}, \; Q_1,$$
where $Q_d$ stands for the set $Q_d=(\Z_{36})_d =\{ x \in L(36) \mid \gcd(x,36)=d \}$. 
Our goal is to describe
$$ G = \bigcap_{k=1}^{7} \aut(\cay(\Z_{36}, B_k)) $$
as the permutation group preserving each of $7$ non-trivial basic Cayley graphs.
It turns out however that we may avoid consideration of all $7$ basic graphs.
(We refer the reader to the texts \cite{FarKM94,KliRRT99,Wie64}
for discussion of corresponding tools, in particular Galois correspondence
between S-rings and permutation groups as well as Schur-Wielandt principle.)

Thus, acting in such spirit, we observe that it is possible to disregard basic sets 
$Q_1, Q_{18}, Q_{12}$, and $Q_9$. Therefore now we define $G$ as group which preserves 
three Cayley graphs
$\Gamma_i$ over $\Z_{36}$ defined by basic sets $Q_2 \cup Q_4$, $Q_6$ and $Q_3$.
These three graphs are conditionally depicted on the three diagrams below
(see also discussion of the rules of the game accepted in these figures).
We admit that ad hoc reasonings are playing a significant role in the ongoing exposition.

\medskip

Graph $\Gamma_1$ is nothing else but a regular graph of valency 12, which has a quotient graph
$\widetilde{\Gamma}_1$ on $12$ metavertices, see Figure~7.
Each metavertex consists of subsets
$\{i,12+i,24+i\}$, where $i \in \Z_{12}$.
Each metaedge substitutes $9$ edges in complete bipartite graphs $K_{3,3}$.
The graphs $\Gamma_1$ and $\widetilde{\Gamma}_1$ have two connectivity components corresponding to
even and odd elements of $\Z_{36}$. An easy way to describe isomorphism type of the components of
$\widetilde{\Gamma}_1$ is $\overline{3 \circ K_2}$, the complement of a $1$-factor on $6$ points.

\bigskip

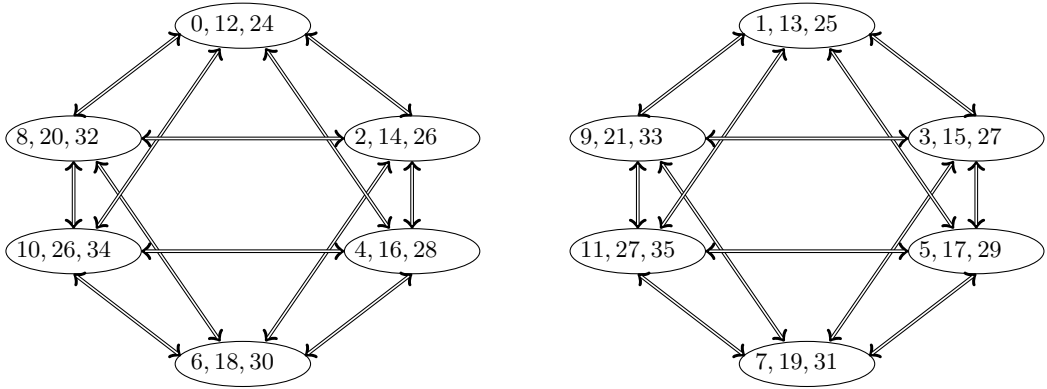
\begin{figure}[h!]
\centering

\begin{tikzpicture}[scale=1.5] ---
\draw (0,0) ellipse (0.6cm and 0.2cm);
\draw (-0.55,0) node[right] {\footnotesize $6,18,30$};
\draw (1.5,1) ellipse (0.6cm and 0.2cm);
\draw (0.9,1) node[right] {\footnotesize $4,16,28$};
\draw (1.5,2) ellipse (0.6cm and 0.2cm);
\draw (0.9,2) node[right] {\footnotesize $2,14,26$};
\draw (0,3) ellipse (0.6cm and 0.2cm);
\draw (-0.55,3) node[right] {\footnotesize $0,12,24$};
\draw (-1.5,2) ellipse (0.6cm and 0.2cm);
\draw (-2.1,2) node[right] {\footnotesize $8,20,32$};
\draw (-1.5,1) ellipse (0.6cm and 0.2cm);
\draw (-2.1,1) node[right] {\footnotesize $10,26,34$};
\draw[<->,double] (0.55,0.07) -- (1.5,0.8);
\draw[<->,double] (-0.55,0.07) -- (-1.5,0.8);
\draw[<->,double] (0.2,0.2) -- (1.3,1.8);
\draw[<->,double] (-0.2,0.2) -- (-1.3,1.8);
\draw[<->,double] (0.55,2.93) -- (1.5,2.2);
\draw[<->,double] (-0.55,2.93) -- (-1.5,2.2);
\draw[<->,double] (0.2,2.8) -- (1.3,1.2);
\draw[<->,double] (-0.2,2.8) -- (-1.3,1.2);
\draw[<->,double] (1.5,1.2) -- (1.5,1.8);
\draw[<->,double] (-1.5,1.2) -- (-1.5,1.8);
\draw[<->,double] (-0.9,1) -- (0.9,1);
\draw[<->,double] (-0.9,2) -- (0.9,2);

\draw (5,0) ellipse (0.6cm and 0.2cm);
\draw (4.45,0) node[right] {\footnotesize $7,19,31$};
\draw (6.5,1) ellipse (0.6cm and 0.2cm);
\draw (5.9,1) node[right] {\footnotesize $5,17,29$};
\draw (6.5,2) ellipse (0.6cm and 0.2cm);
\draw (5.9,2) node[right] {\footnotesize $3,15,27$};
\draw (5,3) ellipse (0.6cm and 0.2cm);
\draw (4.45,3) node[right] {\footnotesize $1,13,25$};
\draw (3.5,2) ellipse (0.6cm and 0.2cm);
\draw (2.9,2) node[right] {\footnotesize $9,21,33$};
\draw (3.5,1) ellipse (0.6cm and 0.2cm);
\draw (2.9,1) node[right] {\footnotesize $11,27,35$};
\draw[<->,double] (5.55,0.07) -- (6.5,0.8);
\draw[<->,double] (4.45,0.07) -- (3.5,0.8);
\draw[<->,double] (5.2,0.2) -- (6.3,1.8);
\draw[<->,double] (4.8,0.2) -- (3.7,1.8);
\draw[<->,double] (5.55,2.93) -- (6.5,2.2);
\draw[<->,double] (4.45,2.93) -- (3.5,2.2);
\draw[<->,double] (5.2,2.8) -- (6.3,1.2);
\draw[<->,double] (4.8,2.8) -- (3.7,1.2);
\draw[<->,double] (6.5,1.2) -- (6.5,1.8);
\draw[<->,double] (3.5,1.2) -- (3.5,1.8);
\draw[<->,double] (4.1,1) -- (5.9,1);
\draw[<->,double] (4.1,2) -- (5.9,2);
\end{tikzpicture}

\caption{$\Gamma_1=\cay(\Z_{36},Q_2 \cup Q_4)$}
\end{figure}

\medskip

Graph $\Gamma_2$ is a disconnected graph of the form $6 \circ C_6$, see Figure~8.
Each cycle $C_6$ is defined on two metavertices from $\Gamma_1$.
Correspondence is observed from diagram.

\bigskip

\begin{figure} [h!]
\centering

\begin{tikzpicture}[scale=1.5] ---
\fill (-0.35,0) circle (1.5pt);
\fill (0,0) circle (1.5pt);
\fill (0.35,0) circle (1.5pt);
\draw (-0.35,-0.1) node[below] {\footnotesize $18$};
\draw (0,-0.1) node[below] {\footnotesize $6$};
\draw (0.35,-0.1) node[below] {\footnotesize $30$};
\fill (1.15,1) circle (1.5pt);
\fill (1.5,1) circle (1.5pt);
\fill (1.85,1) circle (1.5pt);
\draw (1.15,0.9) node[below] {\footnotesize $4$};
\draw (1.5,0.9) node[below] {\footnotesize $16$};
\draw (1.85,0.9) node[below] {\footnotesize $28$};
\fill (1.15,2) circle (1.5pt);
\fill (1.5,2) circle (1.5pt);
\fill (1.85,2) circle (1.5pt);
\draw (1.15,2.1) node[above] {\footnotesize $2$};
\draw (1.5,2.1) node[above] {\footnotesize $14$};
\draw (1.85,2.1) node[above] {\footnotesize $26$};
\fill (-0.35,3) circle (1.5pt);
\fill (0,3) circle (1.5pt);
\fill (0.35,3) circle (1.5pt);
\draw (-0.35,3.1) node[above] {\footnotesize $0$};
\draw (0,3.1) node[above] {\footnotesize $24$};
\draw (0.35,3.1) node[above] {\footnotesize $12$};
\fill (-1.85,2) circle (1.5pt);
\fill (-1.5,2) circle (1.5pt);
\fill (-1.15,2) circle (1.5pt);
\draw (-1.85,2.1) node[above] {\footnotesize $34$};
\draw (-1.5,2.1) node[above] {\footnotesize $10$};
\draw (-1.15,2.1) node[above] {\footnotesize $22$};
\fill (-1.85,1) circle (1.5pt);
\fill (-1.5,1) circle (1.5pt);
\fill (-1.15,1) circle (1.5pt);
\draw (-1.85,0.9) node[below] {\footnotesize $32$};
\draw (-1.5,0.9) node[below] {\footnotesize $8$};
\draw (-1.15,0.9) node[below] {\footnotesize $20$};
\draw (-0.35,3) -- (0,0) -- (0.35,3) -- (-0.35,0) -- (0,3) -- (0.35,0) -- (-0.35,3);
\draw (1.15,2) -- (-1.5,1) -- (1.5,2)  -- (-1.15,1) -- (1.85,2) -- (-1.85,1) -- (1.15,2);
\draw (1.15,1) -- (-1.5,2) -- (1.5,1) -- (-1.15,2) -- (1.85,1) -- (-1.85,2) -- (1.15,1);

\fill (4.65,0) circle (1.5pt);
\fill (5,0) circle (1.5pt);
\fill (5.35,0) circle (1.5pt);
\draw (4.65,-0.1) node[below] {\footnotesize $19$};
\draw (5,-0.1) node[below] {\footnotesize $7$};
\draw (5.35,-0.1) node[below] {\footnotesize $31$};
\fill (6.15,1) circle (1.5pt);
\fill (6.5,1) circle (1.5pt);
\fill (6.85,1) circle (1.5pt);
\draw (6.15,0.9) node[below] {\footnotesize $5$};
\draw (6.5,0.9) node[below] {\footnotesize $29$};
\draw (6.85,0.9) node[below] {\footnotesize $17$};
\fill (6.15,2) circle (1.5pt);
\fill (6.5,2) circle (1.5pt);
\fill (6.85,2) circle (1.5pt);
\draw (6.15,2.1) node[above] {\footnotesize $3$};
\draw (6.5,2.1) node[above] {\footnotesize $15$};
\draw (6.85,2.1) node[above] {\footnotesize $27$};
\fill (4.65,3) circle (1.5pt);
\fill (5,3) circle (1.5pt);
\fill (5.35,3) circle (1.5pt);
\draw (4.65,3.1) node[above] {\footnotesize $1$};
\draw (5,3.1) node[above] {\footnotesize $25$};
\draw (5.35,3.1) node[above] {\footnotesize $13$};
\fill (3.15,2) circle (1.5pt);
\fill (3.5,2) circle (1.5pt);
\fill (3.85,2) circle (1.5pt);
\draw (3.15,2.1) node[above] {\footnotesize $35$};
\draw (3.5,2.1) node[above] {\footnotesize $11$};
\draw (3.85,2.1) node[above] {\footnotesize $23$};
\fill (3.15,1) circle (1.5pt);
\fill (3.5,1) circle (1.5pt);
\fill (3.85,1) circle (1.5pt);
\draw (3.15,0.9) node[below] {\footnotesize $33$};
\draw (3.5,0.9) node[below] {\footnotesize $9$};
\draw (3.85,0.9) node[below] {\footnotesize $21$};
\draw (4.65,3) -- (5,0) -- (5.35,3) -- (4.65,0) -- (5,3) -- (5.35,0) -- (4.65,3);
\draw (6.15,2) -- (3.5,1) -- (6.5,2)  -- (3.85,1) -- (6.85,2) -- (3.15,1) -- (6.15,2);
\draw (6.15,1) -- (3.5,2) -- (6.5,1) -- (3.85,2) -- (6.85,1) -- (3.15,2) -- (6.15,1);
\end{tikzpicture}

\caption{ $\Gamma_2=\cay(\Z_{36},Q_6)$. }
\end{figure}
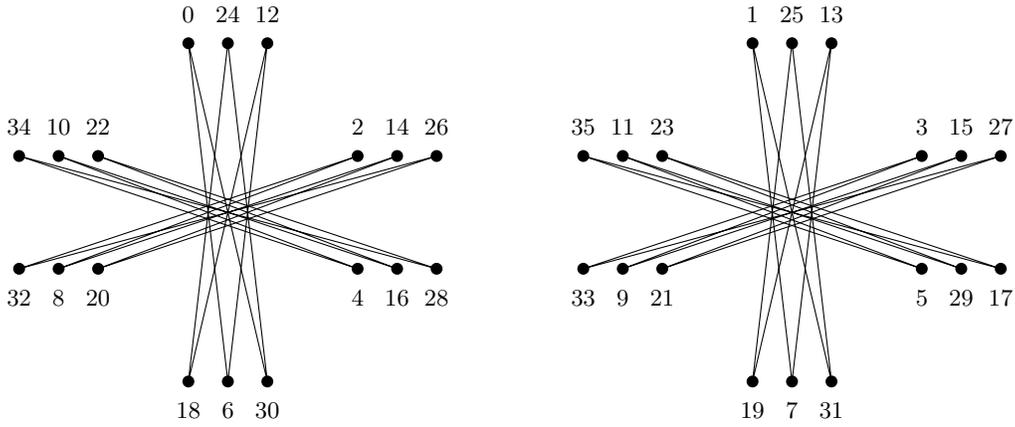

\medskip

Graph $\Gamma_3$ has a more sophisticated nature.
It has three connectivity components defined by the value
of $x \in \Z_{36}$ modulo $3$, one of them is depicted in Figure~9.
Each connectivity component is a bipartite graph with bipartition to odd and even elements.
In addition, each component is $3$-partite with the parts visible on the picture.
Thus finally it may be convenient to regard edge set of a connectivity component as union of edges
from $6$ disjoint quadrangles.

\begin{figure} [h!]
\centering

\begin{tikzpicture}[scale=1.75] ---
\fill (-1.5,0) circle (1.5pt);
\draw (-1.5,0) node[below] {\footnotesize $24$};
\path (-1.5,0) coordinate (24);

\fill (1.5,0) circle (1.5pt);
\draw (1.5,0) node[below] {\footnotesize $30$};
\path (1.5,0) coordinate (30);
\fill (-2,0.5) circle (1.5pt);
\draw (-2,0.5) node[left] {\footnotesize $6$};
\path (-2,0.5) coordinate (6);
\fill (-1,0.5) circle (1.5pt);
\draw (-1,0.5) node[right] {\footnotesize $33$};
\path (-1,0.5) coordinate (33);
\fill (1,0.5) circle (1.5pt);
\draw (1,0.5) node[left] {\footnotesize $21$};
\path (1,0.5) coordinate (21);
\fill (2,0.5) circle (1.5pt);
\draw (2,0.5) node[right] {\footnotesize $12$};
\path (2,0.5) coordinate (12);
\fill (-1.5,1) circle (1.5pt);
\draw (-1.5,1) node[above] {\footnotesize $15$};
\path (-1.5,1) coordinate (15);
\fill (1.5,1) circle (1.5pt);
\draw (1.5,1) node[above] {\footnotesize $3$};
\path (1.5,1) coordinate (3);
\fill (-0.35,2.2) circle (1.5pt);
\draw (-0.35,2.2) node[above] {\footnotesize $0$};
\path (-0.35,2.2) coordinate (0);
\fill (-0.35,1.5) circle (1.5pt);
\draw (-0.35,1.5) node[below] {\footnotesize $9$};
\path (-0.35,1.5) coordinate (9);
\fill (0.35,2.2) circle (1.5pt);
\draw (0.35,2.2) node[above] {\footnotesize $18$};
\path (0.35,2.2) coordinate (18);
\fill (0.35,1.5) circle (1.5pt);
\draw (0.35,1.5) node[below] {\footnotesize $27$};
\path (0.35,1.5) coordinate (27);

\draw (0) -- (3) (0) -- (21) (18) -- (3) (18) -- (21);
\draw (0) -- (15) (0) -- (33);
\draw (9) -- (12) (9) -- (24) (9) -- (30) (9) -- (6);
\draw (18) -- (33) (18) -- (15);
\draw (27) -- (30) (27) -- (6) (27) -- (12) (27) -- (24);
\draw (3) -- (6) (3) -- (24);
\draw (12) -- (15) (12) -- (33);
\draw (21) -- (24) (21) -- (6);
\draw (30) -- (33) (30) -- (15);

\end{tikzpicture}

\caption{ A connectivity component of $\Gamma_3=\cay(\Z_{36},Q_3)$.}
\end{figure}
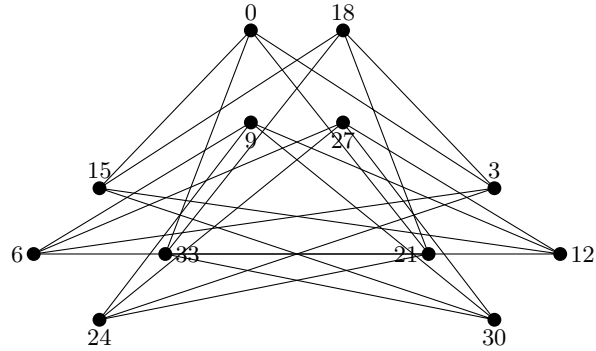

\medskip

Now we are prepared to claim that the desired group $G$ has the following structure: 
$$ G = \Z_{2}^{6} \centerdot \big( \, (S_3 \wr S_3) \centerdot \Z_2 \, \big),$$
and thus it has order $2^6 \cdot 6^4 \cdot 2 = 2^{11} \cdot 3^4.$
To prove this claim, we will present concrete automorphisms from $G$, will comment their action on
the basic graphs, and will count the order of the group, generated by these permutations.
First we wish to describe $64$ permutations from $G$, which preserve each metavertex of $\Gamma_1$
and each connectivity component of $\Gamma_2$.
(Of course, in addition, they preserve the remaining graph, this time $\Gamma_3$.)
In fact, we restrict ourselves by list of $3$ permutations which are corresponding to the
connectivity component of $\Gamma_3$ given in Figure~8.
\begin{eqnarray*}
g^{(1)}_1 &=& (3,15)(6,30)(12,24)(21,33), \\
g^{(1)}_2 &=& (0,12)(18,30)(9,21)(3,27), \\
g^{(1)}_3 &=& (0,24)(6,18)(9,33)(15,27).
\end{eqnarray*}
Similarly, two more sets of permutations $g^{(i)}_1,g^{(i)}_2,g^{(i)}_3, i=2,3,$ are defined.
Altogether, involutions from three groups, isomorphic to $(\Z_2)^2$ are listed.  
Direct product of these three groups provides group $(\Z_2)^6,$ 
forming first factor in description of $G$.

Now we wish to justify part of the formula $S_3 \wr S_3$. It is helpful to think about the group acting faithfully on the set of $9$ anti-cliques of size $4$ visible from the diagram of $\Gamma_3$.
First, consider permutations on $\Z_{36}$ defined as
\begin{eqnarray*}
g_4 & \colon & x \mapsto x+4, \text{ and } \\
g_5 &=& (0)(1,35)(2,34)(3,33) \cdots (17,19)(18).
\end{eqnarray*}
Clearly, these permutations generate a subgroup, which acts as $S_3$ on the connected components of
$\Gamma_3$ and preserves odd and even parts.
On next step, consider
\begin{eqnarray*}
g^{(1)}_6 &=& (0,6,12,18,24,30)(3,9,15,21,27,33), \text{ and } \\ 
g^{(1)}_7  &=& g^{(1)}_1 = (3,15)(6,30)(12,24)(21,33). 
\end{eqnarray*}
Check that $\sg{g^{(1)}_6, g^{(1)}_7}$ acts as $S_3$ on 
the connected component of $\Gamma_3$ given in Figure~8,
it preserves other components, and of course it is an automorphism group of two remaining basic graphs.
Similarly, two more sets of permutations $\sg{g^{(i)}_6,g^{(i)}_7}, i=2,3,$ are defined.
Last natural permutation on $\Z_{36}$ is defined as $g_8 \colon x \mapsto x + 1$, which
clearly interchanges odd and even vertices.

We suggest the reader to check that the permutations $g^{(i)}_1,g^{(i)}_2, \ldots, g_8$ exposed above
(which belong to $G$ indeed) generate the group of the desired order $2^{11} \cdot 3^{4}$.
It is a standard exercise in computational algebraic graph theory to confirm that we already
encountered the entire group $G$. \QED
\end{exm}

In next sections group $G$ will appear again, though in different incarnations, thus
helping the reader again and again to build a bridge between our theoretical reasonings and practical ad hoc computations.

\section{Crested products}

In this section we provide a short digest of the paper \cite{BaiC05}, which is adopted essentially
for the purposes of the current presentation.  We refer to \cite{BaiC05} for accurate proofs of the claims
presented below, while ongoing level of rigor follows the intuitive style of the previous section.

Recall that our foremost goal is to investigate and to extend the possibility to build arbitrary sublattice
$L$ of $L(n)$ from trivial lattices using only simple reduction rules.
The {\em trivial sublattice} of $L(n)$ consists of only the elements
$1$ and $n$, and it will be denoted by $T_n$. We may  prove that such ``easy life''
is possible if and only if $n=p^e$ or $n=p^eq$ or $n=pqr$ for distinct primes $p,q$ and $r$ (see Section 11).
We wish to define binary operation $\otimes_d, d\in \N,$ for lattices with the following goals in mind.

\begin{itemize}
\item Special cases of $\otimes_d$ give back simple reduction rules.
\item Every sublattice $L$ of $L(n)$ such that $1,n \in L$ can be built from trivial lattices using only operations $\otimes_d$.
\item If $L$ is built from trivial lattices as
$$ L = T_{d_k} \otimes_{d_{k-1}}
\Big( T_{d_{k-1}} \otimes_{d_{k-2}} ( \cdots   \otimes_{d_2} (T_{d_2} \otimes_{d_1} T_{d_1}) \cdots )\Big),$$
then $\aut(L)$ can be nicely
described in terms of symmetric groups $\aut(T_{d_i})$ $=S_{d_i}$.
\end{itemize}

\noindent In what follows this desired operation $\otimes_d$ will be called {\em crested product}.
The word ``crested'', suggested in \cite{BaiC05}, is a mixture of ``crossed'' and ``nested'',
and is also cognate with the meaning of ``wreath'' in ``wreath product''.
Due to the existence of the
bijections between S-rings of traces over $\Z_n$, sublattices of $L(n)$,
rational association schemes (invariant with respect to regular cyclic groups)
and orthogonal group block structures on $\Z_n$,
the desired new operation may be translated in a few corresponding diverse languages. We
prefer to start with orthogonal block structures (see \cite[Definition~3]{BaiC05}).

\medskip

For $i=1,2$, let $F_i$ be a partition of $X_i$. Define $F_1 \times F_2$ to be
the partition of $X_1 \times X_2$ whose adjacency matrix $A(F_1 \times F_2)$ is 
$A(F_1) \otimes A(F_2)$ (see 2.4).

\begin{defi}\label{def:crestedp}
For $i=1,2$, let $\F_i$ be an orthogonal block structure on a set $X_i$, and let $F_i \in \F_i$.
The {\em (simple) crested product} of $\F_1$ and $\F_2$ with respect to $F_1$ and $F_2$ is the following set
$\P$ of partitions of $X_1 \times X_2$:
$$ \P = \big\{ \, P_1 \times P_2 \, \mid \, P_1 \in \F_1, P_2 \in \F_2, P_1 \sqsubseteq F_1 \text{ or }
P_2 \sqsupseteq F_2 \, \big\}. $$
\end{defi}

\medskip

It can be proved that the crested product, as just defined, is an orthogonal block structure.
The reader may be easily convinced that indeed, crossing and nesting are special cases of the
crested product. An important subclass of orthogonal block structures consists of the
poset block structures. It can be proved that crested products of poset block structures remain poset
block structures. Moreover, every poset block structure can be attained from trivial block structures by a
repeated use of crested products. Thus it can be proved that the crested products satisfy the above
three goals. (Note that our claim about the fulfillment of the above goals literally is actual for the poset block structures on $\Z_n$. We avoid discussion of difficulties, which may appear in more general cases.)

\medskip

The formal definition of crested product $\otimes_{d}$ (adopted for the orthogonal 
group block structures on $\Z_n$) is as follows.

\begin{defi}\label{def:crestedp2}
For $i=1,2$, let $n_i \in \N$, $L_i$ be a sublattice of $L(n_i)$ such that $1,n_i \in L(n_i)$,
and $d$ be in $L_2$ such that $\gcd(n_1,n_2/d)=1$. Then the sublattice
$L_1 \otimes_d L_2$ of $L(n_1n_2)$ is defined as
$$ L_1 \otimes_d L_2 = \big\{ \,  l_1l_2 \, \mid \, l_1=1, l_2\in L_2, \text{ or } l_1 \in L_1,
l_2 \in L_2 \text{ with } d \mid l_2 \, \big\}. $$
\end{defi}

\medskip

\noindent Notice that, operations $\otimes_d$ include simple reduction rules 1 and 2 as special cases.
Namely, in case $d=n_2,$ and $L_1=T_{n_1}$ or $L_2=T_{n_2}$ we get reduction rule 1,
and in case $d=1$ reduction rule 2.

Consider the orthogonal group block structure on $\Z_{n_1n_2}$
corresponding to the lattice $L_1 \otimes_d L_2$. This is weakly isomorphic to the
crested product of the block structure on $\Z_{n_1}$ corresponding to $L_1$
and that one on $\Z_{n_2}$ corresponding to $L_2$ with respect to partitions $F_{Z_1}$ and $F_{Z_d}$
in the sense of Definition~\ref{def:crestedp}, justifying the name ``crested product'' for
$\otimes_d$.

\begin{exm}\label{ex:n=36'} ({\it Example~\ref{ex:n=36} revised}.)
Let $L$ be the sublattice of $L(36)$ given in Figure~6.
To each of $8$ nodes in diagram for $L$ naturally corresponds a partition of
$\Z_{36}$.  Because $L$ is a lattice, we get a corresponding orthogonal block structure on $\Z_{36}$.
Naive description of nodes of $L$ looks as follows: consider all nodes in $L$ and take into consideration
those ones which are in $L\subl{18}$ or are larger than $d=2$. 
Let $L_1=\{1,2\}$ on $\Z_2$ and $L_2=L\subl{18}=\{1,2,3,6,18\}$ on $\Z_{18}$.
Then by Definition~\ref{def:crestedp2} we obtain
$$ L = \{1 \cdot 1,1 \cdot 2, 1 \cdot 3, 1 \cdot 6, 1 \cdot 18, 2 \cdot 2, 2 \cdot 6, 2 \cdot
18\} = L_1 \otimes_2 L_2. $$
Moreover, using properly notation for the crested product of lattices, the product 
$L_1 \otimes_2 L_2$ is depicted in Figure~10.

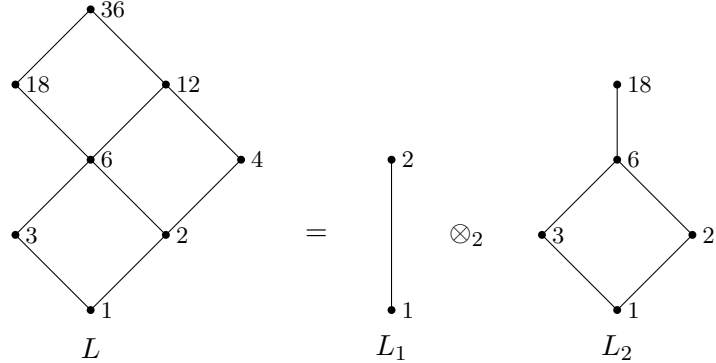
\begin{figure}[h!]
\centering

\begin{tikzpicture}[scale=1] ---
\fill (0,0) circle (1.5pt);
\draw (0,0) node[right] {\footnotesize $1$};
\fill (-1,1) circle (1.5pt);
\draw (-1,1) node[right] {\footnotesize $3$};
\fill (1,1) circle (1.5pt);
\draw (1,1) node[right] {\footnotesize $2$};
\fill (0,2) circle (1.5pt);
\draw (0,2) node[right] {\footnotesize $6$};
\fill (2,2) circle (1.5pt);
\draw (2,2) node[right] {\footnotesize $4$};
\fill (1,3) circle (1.5pt);
\draw (1,3) node[right] {\footnotesize $12$};
\fill (-1,3) circle (1.5pt);
\draw (-1,3) node[right] {\footnotesize $18$};
\fill (0,4) circle (1.5pt);
\draw (0,4) node[right] {\footnotesize $36$};
\draw (0,0) -- (-1,1) -- (0,2) -- (1,1) -- (0,0);
\draw (1,1) -- (2,2) -- (1,3) -- (0,2);
\draw (1,3) -- (0,4) -- (-1,3) -- (0,2);
\draw (0,-0.5) node {$L$};

\draw (3,1) node {$=$};

\fill (4,0) circle (1.5pt);
\draw (4,0) node[right] {\footnotesize $1$};
\fill (4,2) circle (1.5pt);
\draw (4,2) node[right] {\footnotesize $2$};
\draw (4,0) -- (4,2);
\draw (4,-0.5) node {$L_1$};

\draw (5,1) node {$\otimes_2$};

\fill (7,0) circle (1.5pt);
\draw (7,0) node[right] {\footnotesize $1$};
\fill (8,1) circle (1.5pt);
\draw (8,1) node[right] {\footnotesize $2$};
\fill (7,2) circle (1.5pt);
\draw (7,2) node[right] {\footnotesize $6$};
\fill (6,1) circle (1.5pt);
\draw (6,1) node[right] {\footnotesize $3$};
\fill (7,3) circle (1.5pt);
\draw (7,3) node[right] {\footnotesize $18$};
\draw (7,0) -- (8,1) -- (7,2) -- (6,1) -- (7,0) (7,2) -- (7,3);
\draw (7,-0.5) node {$L_2$};
\end{tikzpicture}

\caption{Decomposition $L = L_1 \otimes_2 L_2$.}
\end{figure}

We now easily interpret $L$ with the aid of Definition~\ref{def:crestedp} as crested product.
Namely, consider subgroups $Z_m \le \Z_{36}$ for $m=2,4$ and $18$.
We have $Z_2=\{0,18\},$ $Z_4=\{0,9,18,27\},$ and write the 
quotient group $Z_4/Z_2$ as $Z_4/Z_2=\{Z_2,Z_2+9\}$. As $0$ and $9$ form a complete set of 
coset representatives of the subgroup $Z_{18}$ in $\Z_{36},$ every element $x$ in $\Z_{36}$ can be written
uniquely as
$$ x = x_1 + x_2, \text{ where } x_1 \in \{0,9\} \text{ and } x_2 \in Z_{18}.$$
Therefore, we can define the bijective mapping
$$ f \colon \Z_{36} \to Z_4/Z_2 \times Z_{18}, \, x \mapsto (Z_2+x_1,x_2).$$
Let $F_1$ be the trivial partition $\{\, \{Z_2\},\{Z_2+9\}\, \}$ of $Z_4/Z_2,$ and
$F_2$ be the partition isomorphic to $9 \circ K_2$ of $Z_{18}$ into cosets of $Z_2$.
Here and later on we freely identify partitions with the graphs defined by the corresponding
equivalence relations. We leave for the reader to check that $f$ is a week isomorphism from
our block structure $L=L_1 \otimes_2 L_2$
to the crested product of the block structure on $Z_4/Z_2$ corresponding to $L_1,$ and the block structure on
$Z_{18}$ corresponding to $L_2$ with respect to $F_1$ and $F_2$.

Eventually, notice that simple reduction rules apply to $L_2$. We obtain that 
$L_2 = T_3 \otimes_6 (T_2 \otimes_1 T_3),$ therefore, $L$ actually decomposes as
$$ 
L = T_2 \otimes_2 (T_3 \otimes_6 (T_2 \otimes_1 T_3)).
$$ \QED
\end{exm}

\medskip

In the rest of the section we turn to the group $\aut(L_1 \otimes_d L_2)$.
It remains to translate everything to the language of association schemes, and after that the one
of permutation groups, with the goal that finally $\aut(L_1 \otimes_d L_2)$ is described in terms of
$\aut(L_1)$ and $\aut(L_2)$.  We refer again to the paper \cite{BaiC05}, where such goal is fulfilled to a
certain extent.  Namely, it is proved that for the case of poset block structures one gets that crested product of  $\aut(L_1)$ and $\aut(L_2)$ preserves the crested product of $L_1$ and $L_2$.
Instead of a discussion of corresponding precise definitions and formulations, we prefer to play again on
the level of our striking example.

\begin{exm}\label{ex:n=36''} ({\it Continuation of Example~\ref{ex:n=36'}}.)
We again use freely the possibility to switch at any moment between languages of lattices, S-rings, and association schemes.  In the above notation we get $L_1 = \{1,2\}$ on $\Z_2$, $L_2 = \{1,2,3,6,18\}$ on $\Z_{18}\}$ and
$L = \{1 \cdot 1,1 \cdot 2, 1 \cdot 3, 1 \cdot 6, 1 \cdot 18, 2 \cdot 1, 2 \cdot 2, 2 \cdot 6, 2 \cdot
18\}$ on $\Z_{36}$.

In our previous attempt it was natural and convenient to consider automorphism groups of basic graphs
(regarded as rational circulant graphs).  We proceeded finally with three such graphs.  At the current stage we see
$G = \aut(L)$ with the aid of group basis in the corresponding S-ring (cf. \cite{BriM79,Muz93}).
Clearly, each element of a group basis corresponds to a partition of $\Z_{36}$ into cosets of a suitable subgroup.  Therefore, now we get
$$ G = \bigcap_{l \in L} \aut(\cay(\Z_{36}, Z_m)),$$
where $Z_m$ is the unique subgroup of $\Z_{36}$ of order $m$.
Thus we have immediately that in fact $G$ is the automorphism group of four partitions defined by
$Z_m,$ namely $m=2,3,4$ and $18$. We again describe this group, using a suitable diagram below which exhibits
simultaneously all the partitions.

\begin{figure} 
\centering

\begin{tikzpicture}[scale=1.5] ---
\fill (-1,0) circle (1.5pt);
\draw (-1,0) node[left] {\footnotesize $26$};
\path (-1,0) coordinate (26);
\fill (-1,0.5) circle (1.5pt);
\draw (-1,0.5) node[left] {\footnotesize $14$};
\path (-1,0.5) coordinate (14);
\fill (-1,1) circle (1.5pt);
\draw (-1,1) node[left] {\footnotesize $2$};
\path (-1,1) coordinate (2);
\fill (-0.5,0) circle (1.5pt);
\draw (-0.5,0) node[right] {\footnotesize $8$};
\path (-0.5,0) coordinate (8);
\fill (-0.5,0.5) circle (1.5pt);
\draw (-0.5,0.5) node[right] {\footnotesize $32$};
\path (-0.5,0.5) coordinate (32);
\fill (-0.5,1) circle (1.5pt);
\draw (-0.5,1) node[right] {\footnotesize $20$};
\path (-0.5,1) coordinate (20);
\draw (8) -- (26) -- (2) -- (20) -- (8) (14) -- (32);

\fill (1,0) circle (1.5pt);
\draw (1,0) node[right] {\footnotesize $17$};
\path (1,0) coordinate (17);
\fill (1,0.5) circle (1.5pt);
\draw (1,0.5) node[right] {\footnotesize $5$};
\path (1,0.5) coordinate (5);
\fill (1,1) circle (1.5pt);
\draw (1,1) node[right] {\footnotesize $29$};
\path (1,1) coordinate (29);
\fill (0.5,0) circle (1.5pt);
\draw (0.5,0) node[left] {\footnotesize $35$};
\path (0.5,0) coordinate (35);
\fill (0.5,0.5) circle (1.5pt);
\draw (0.5,0.5) node[left] {\footnotesize $23$};
\path (0.5,0.5) coordinate (23);
\fill (0.5,1) circle (1.5pt);
\draw (0.5,1) node[left] {\footnotesize $11$};
\path (0.5,1) coordinate (11);
\draw (35) -- (17) -- (29) -- (11) -- (35) (23) -- (5);

\fill (-1,1.5) circle (1.5pt);
\draw (-1,1.5) node[left] {\footnotesize $34$};
\path (-1,1.5) coordinate (34);
\fill (-1,2) circle (1.5pt);
\draw (-1,2) node[left] {\footnotesize $22$};
\path (-1,2) coordinate (22);
\fill (-1,2.5) circle (1.5pt);
\draw (-1,2.5) node[left] {\footnotesize $10$};
\path (-1,2.5) coordinate (10);
\fill (-0.5,1.5) circle (1.5pt);
\draw (-0.5,1.5) node[right] {\footnotesize $16$};
\path (-0.5,1.5) coordinate (16);
\fill (-0.5,2) circle (1.5pt);
\draw (-0.5,2) node[right] {\footnotesize $4$};
\path (-0.5,2) coordinate (4);
\fill (-0.5,2.5) circle (1.5pt);
\draw (-0.5,2.5) node[right] {\footnotesize $28$};
\path (-0.5,2.5) coordinate (28);
\draw (16) -- (34) -- (10) -- (28) -- (16) (4) -- (22);

\fill (1,1.5) circle (1.75pt);
\draw (1,1.5) node[right] {\footnotesize $7$};
\path (1,1.5) coordinate (7);
\fill (1,2) circle (1.5pt);
\draw (1,2) node[right] {\footnotesize $31$};
\path (1,2) coordinate (31);
\fill (1,2.5) circle (1.5pt);
\draw (1,2.5) node[right] {\footnotesize $19$};
\path (1,2.5) coordinate (19);
\fill (0.5,1.5) circle (1.5pt);
\draw (0.5,1.5) node[left] {\footnotesize $25$};
\path (0.5,1.5) coordinate (25);
\fill (0.5,2) circle (1.5pt);
\draw (0.5,2) node[left] {\footnotesize $13$};
\path (0.5,2) coordinate (13);
\fill (0.5,2.5) circle (1.5pt);
\draw (0.5,2.5) node[left] {\footnotesize $1$};
\path (0.5,2.5) coordinate (1);
\draw (25) -- (7) -- (19) -- (1) -- (25) (13) -- (31);

\fill (-1,3) circle (1.5pt);
\draw (-1,3) node[left] {\footnotesize $24$};
\path (-1,3) coordinate (24);
\fill (-1,3.5) circle (1.5pt);
\draw (-1,3.5) node[left] {\footnotesize $12$};
\path (-1,3.5) coordinate (12);
\fill (-1,4) circle (1.5pt);
\draw (-1,4) node[left] {\footnotesize $0$};
\path (-1,4) coordinate (0);
\fill (-0.5,3) circle (1.5pt);
\draw (-0.5,3) node[right] {\footnotesize $6$};
\path (-0.5,3) coordinate (6);
\fill (-0.5,3.5) circle (1.5pt);
\draw (-0.5,3.5) node[right] {\footnotesize $30$};
\path (-0.5,3.5) coordinate (30);
\fill (-0.5,4) circle (1.5pt);
\draw (-0.5,4) node[right] {\footnotesize $18$};
\path (-0.5,4) coordinate (18);
\draw (6) -- (24) -- (0) -- (18) -- (6) (12) -- (30);

\fill (1,3) circle (1.5pt);
\draw (1,3) node[right] {\footnotesize $15$};
\path (1,3) coordinate (15);
\fill (1,3.5) circle (1.5pt);
\draw (1,3.5) node[right] {\footnotesize $3$};
\path (1,3.5) coordinate (3);
\fill (1,4) circle (1.5pt);
\draw (1,4) node[right] {\footnotesize $27$};
\path (1,4) coordinate (27);
\fill (0.5,3) circle (1.5pt);
\draw (0.5,3) node[left] {\footnotesize $33$};
\path (0.5,3) coordinate (33);
\fill (0.5,3.5) circle (1.5pt);
\draw (0.5,3.5) node[left] {\footnotesize $21$};
\path (0.5,3.5) coordinate (21);
\fill (0.5,4) circle (1.5pt);
\draw (0.5,4) node[left] {\footnotesize $9$};
\path (0.5,4) coordinate (9);
\draw (33) -- (15) -- (27) -- (9) -- (33) (21) -- (3);

\draw[gray] (0,-0.25) -- (0,4.25) (18) -- (9) (30) -- (21) (6) -- (33) (28) -- (1) (4) -- (13) (16) -- (25)
(20) -- (11) (32) -- (23) (8) -- (35);

\end{tikzpicture}

\caption{Coset-partitions of $\Z_{36}$ defined by $Z_m$ for $m=2,3,4,18$.}
\end{figure}
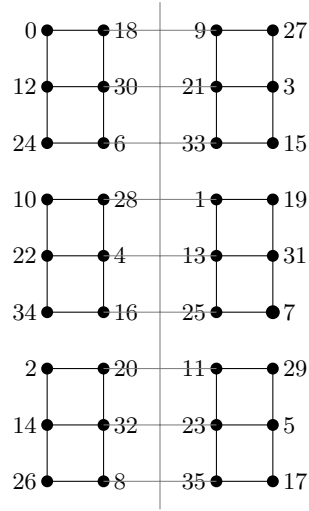

Comments about the diagram. First partition $2 \circ K_{18}$ (due to $Z_{18}$) is presented by division to
left and right part (odd and even numbers).  Horizontal lines represent $9 \circ K_4$ (due to $Z_4$).
Finally, we have $6$ connected components of size $6$.  Columns of all such components entirely
provide $12 \circ K_3$ (due to $Z_3$), while rows give $18 \circ K_2$ (due to $Z_2$).

We now describe the automorphism group $G$ as follows
$G = \widetilde{G} \centerdot S_2$, where $\widetilde{G}$ is the stabilizer of left part of
the picture (clearly left and right parts may be exchanged).
Left part is a wreath product of the groups of the three components.  Thus we get
$G = (S_3 \wr \widehat{G}) \centerdot S_2,$ where $\widehat{G}$ is the stabilizer of a
component. Stabilizer of left upper component, according to simple rule, is $S_2 \times S_3$, and in addition,
an independent copy of $S_2$ transposes columns in corresponding right part of the upper component.
We have thus obtained the formula
$$ G=(S_3 \wr ((S_3 \times S_2) \times S_2)) \centerdot S_2, $$
with the order $|G| = 2 \cdot 3! \cdot 24^3 = 2^{11} \cdot 3^4$.
We expect that the reader will admit that the current arguments are more transparent and straightforward,
however, we again are depending on the use of ad hoc tricks of geometrical and combinatorial nature.

\medskip

It turns out that the above argumentation may be modified into certain nice formal rule with the aid of the use of crested product, taking into account the decomposition formula presented for the lattice
$L$ in the consideration, that is $L = L_1 \otimes_2 L_2$.

Regarding as sets, let $\Z_{36} = \Z_2 \times \Z_{18}$.
In definition of crested product first ingredient corresponds to active while second to passive groups.
Thus in our case $G$ is regarded as a subgroup of the wreath product
$\aut(L_1) \wr \aut(L_2)=G_1 \wr G_2$ or more precisely as
$B \rtimes G_1$, where $B$ is base group and $G_1$ is top group.
Note that at this stage $B$ is just a subgroup of the base group, corresponding to the usual wreath product.  
Using our toolkit of simple rules, we obtain that
$$ 
G_1=\aut(L_1)=S_2, \text{ and }
G_2=\aut(L_2)=S_3 \wr (S_2 \times S_3).
$$ 
We have to understand the structure of the base group.  Recall that in our case $B$ is subgroup of group
$G_2^{\Z_2}$.  To describe $B$ we refer to the partition $F_2=F_{Z_2}$ of $\Z_{18}$, which is preserved by
$G_2$.  Clearly, this partition $F_2$ is of the kind $9 \circ K_2$.
Note that we have also a trivial partition $F_1=F_{Z_1}$ of kind $2 \circ K_1$ which is preserved by $G_1$.
Now we are looking for the subgroup $N$ of group $G_2$ which fixes each part of the partition $F_2$.
Clearly, in our case $N$ is isomorphic to $S_2^3$.

It turns out (see \cite{BaiC05} for general justification) that the base group $B$ is generated by
$N^{F_1}$ and $G_2$.  Here $N^{F_1}$ is embedded in $G_{2}^{\Z_2}$ as the set of functions which are constant on the classes of $F_1$ and take values diagonally.  $G_2$ is embedded diagonally.  $G_2$ normalizes $N^{F_1}$, therefore
their product is a group,  while intersection is $N$.  Thus we obtain that
$$ |B| = \frac{ |N^{F_1}| \cdot |G_2|}{|N|}= \frac{|N|^2 \cdot G_2}{|N|} = |N| \cdot |G_2|. $$
Finally we get the order of the group $G=B \rtimes G_1$ as
$|G_1| \cdot |N| \cdot |G_2| = 2 \cdot 8 \cdot \, 3!(2 \cdot 3 )^3 = 2^{11} \cdot 3^4$, as desired.
(In fact, again we first obtain that the automorphism group of $L$ has order at least $2^{11} \cdot 3^4 $.
After that, exactly like in Example~\ref{ex:n=36}, we have to check that
$B \rtimes G_1$ indeed coincides with the entire group $G$.) \QED
\end{exm}

\noindent{\bf Remark.}
We wish to use an extra chance to explain the role of index 2 in our notation for the used version of crested product.  Hopefully, the following pictorial explanation (see Figure~12) may help. Here selected node in $L_2$ is origin of 
the index. We multiply part of $L_1$ strictly below the index on $L_2$, after that $L_1$ on the part of $L_2$ 
above the index and amalgamate the two products.

\begin{figure}[h!]
\centering

\begin{tikzpicture}[scale=0.75] ---

\fill (0,0.5) circle (1.5pt);
\draw (0,0.5) node[right] {\footnotesize $1$};
\fill[gray] (0,1.5) circle (2.5pt);
\draw (0,1.5) node[right] {\footnotesize $2$};
\draw (0,0.5) -- (0,1.5);

\draw (1,1) node {$\otimes_2$};

\fill (3,0) circle (1.5pt);
\draw (3,0) node[right] {\footnotesize $1$};
\fill[gray] (4,1) circle (2.5pt);
\draw (4,1) node[right] {\footnotesize $2$};
\fill (3,2) circle (1.5pt);
\draw (3,2) node[right] {\footnotesize $6$};
\fill (2,1) circle (1.5pt);
\draw (2,1) node[right] {\footnotesize $3$};
\fill (3,3) circle (1.5pt);
\draw (3,3) node[right] {\footnotesize $18$};
\draw (3,0) -- (4,1) -- (3,2) -- (2,1) -- (3,0) (3,2) -- (3,3);

\draw (5,1) node {$=$};

\fill (7,0) circle (1.5pt);
\draw (7,0) node[right] {\footnotesize $1$};
\fill (8,1) circle (1.5pt);
\draw (8,1) node[right] {\footnotesize $2$};
\fill (7,2) circle (1.5pt);
\draw (7,2) node[right] {\footnotesize $6$};
\fill (6,1) circle (1.5pt);
\draw (6,1) node[right] {\footnotesize $3$};
\fill (7,3) circle (1.5pt);
\draw (7,3) node[right] {\footnotesize $18$};
\draw (7,0) -- (8,1) (7,2) -- (6,1) -- (7,0);
\draw[dotted] (8,1) -- (7,2) -- (7,3);

\draw (9.25,1) node {$\cup$};

\fill (9.5,3) circle (1.5pt);
\draw (9.5,3) node[right] {\footnotesize $18$};
\fill (10.5,2) circle (1.5pt);
\draw (10.5,2) node[right] {\footnotesize $6$};
\fill (10.5,4) circle (1.5pt);
\draw (10.5,4) node[right] {\footnotesize $36$};
\fill (11.5,1) circle (1.5pt);
\draw (11.5,1) node[right] {\footnotesize $2$};
\fill (11.5,3) circle (1.5pt);
\draw (11.5,3) node[right] {\footnotesize $12$};
\fill (12.5,2) circle (1.5pt);
\draw (12.5,2) node[right] {\footnotesize $4$};
\draw (11.5,1) -- (12.5,2) -- (11.5,3) -- (10.5,4) -- (9.5,3) (10.5,2) -- (11.5,3);
\draw[dotted] (11.5,1) -- (10.5,2) -- (9.5,3);

\draw (5,-3) node {$=$};

\fill (7,-5) circle (1.5pt);
\draw (7,-5) node[right] {\footnotesize $1$};
\fill (6,-4) circle (1.5pt);
\draw (6,-4) node[right] {\footnotesize $3$};
\fill (8,-4) circle (1.5pt);
\draw (8,-4) node[right] {\footnotesize $2$};
\fill (7,-3) circle (1.5pt);
\draw (7,-3) node[right] {\footnotesize $6$};
\fill (9,-3) circle (1.5pt);
\draw (9,-3) node[right] {\footnotesize $4$};
\fill (8,-2) circle (1.5pt);
\draw (8,-2) node[right] {\footnotesize $12$};
\fill (6,-2) circle (1.5pt);
\draw (6,-2) node[right] {\footnotesize $18$};
\fill (7,-1) circle (1.5pt);
\draw (7,-1) node[right] {\footnotesize $36$};

\draw (7,-5) -- (6,-4) -- (7,-3) -- (8,-4) -- (7,-5);
\draw (8,-4) -- (9,-3) -- (8,-2) -- (7,-1) -- (6,-2) -- (7,-3) -- (8,-2); 
\end{tikzpicture}

\caption{``Amalgamation of lattices"}
\end{figure}
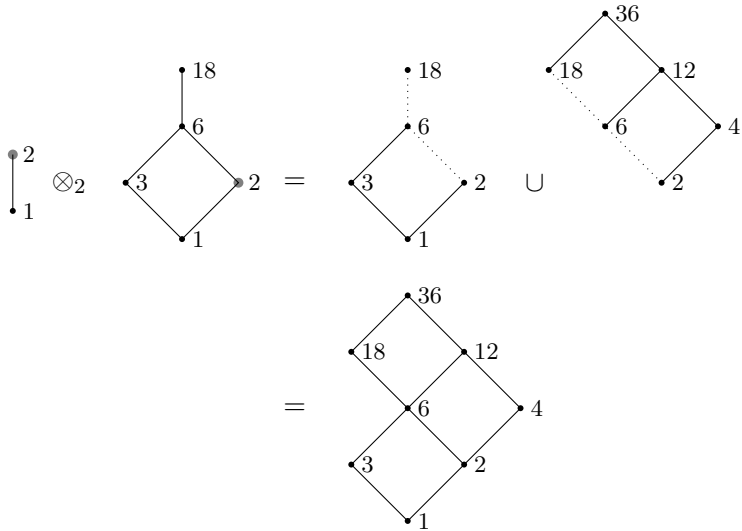


\medskip

In our eyes the formulated goal to create for the reader a context with the aid of an example is fulfilled.
In principle, based on the earned experience, one can go ahead and prove that the leads formulated above are
completely fulfilled with the aid of the crested product.
However, this will not be done in the current paper. We finish the deviation, developed in Sections~7-8,
and return to the main stream of the presentation, aiming to exploit another classical generalization of
the operation of wreath product.

\section{Generalized wreath products}

Let $(I,\preceq)$ be a poset. A subset $J \subseteq I$ is called \emph{ancestral}
if $i \in J$ and $i \preceq j$ imply that $j \in J$ for all $i,j \in I$.
For $i \in I$, we put $A_i$ for the ancestral subset
$$ A_i=\{j \in I \, \mid \, i \prec j\}.$$
The set of all ancestral subsets of $I$ we denote by $\anc((I,\preceq))$.
For each $i \in I$, fix a set $X_i$ of cardinality at least $2$, and let $X=\prod_{i \in I}X_i$.
We write elements $x$ in $X$ as $x=(x_i)_{i \in I}$ or simply as $x=(x_i)$.
For $J \subseteq I$, let $\sim_J$ be the equivalence relation on $X$ given as
$$ (x_i) \sim_J (y_i) \iff x_j=y_j \text{ for all } j \in J,$$
and denote by $\Pi(J)$  the corresponding partition of $X$.
Now, the {\em poset block structure} defined by the poset $(I,\preceq)$ and the sets $X_i$ is
the block structure on $X$ consisting of all partitions $\Pi(J)$ that $J \in \anc((I,\preceq))$.
Denote $\F$ this block structure.
Let $J,J' \in \anc((I,\preceq))$. Both sets $J \cap J'$ and $J \cup J'$ are ancestral, and we have
$$ \Pi(J) \wedge \Pi(J')=\Pi(J \cup J') \text{ and }
\Pi(J) \vee \Pi(J')=\Pi(J \cap J'). $$
Thus the poset block structure $\F$ is an orthogonal block structure.
Further, the equivalence $\Pi(J) \sqsubseteq \Pi(J') \iff J' \subseteq J$ holds, and the
mapping $J \mapsto \Pi(J)$ is an anti-isomorphism from the
lattice $(\anc((I,\preceq)),\subseteq)$ to the lattice $(\F,\sqsubseteq)$ (thus these have Hasse
diagrams dual to each other).  The lattice $(\anc((I,\preceq)),\subseteq)$ is obviously
distributive, and by the previous remarks so is $\F$. The following converse is due to
Bailey and Speed \cite{SpeB82} (see also \cite[Theorem 5]{Bai81}).

\begin{thm}\label{thm:BS}
An orthogonal block structure is distributive if and only if it is weakly isomorphic to
a poset block structure.
\end{thm}

Note that, in particular, the orthogonal group block structures on $\Z_n$ are poset block
structures. We continue consideration of our striking example.

\begin{exm}\label{exm:n=36'''} ({\it Example~\ref{ex:n=36} revised}.)
Let $L$ be the sublattice of $L(36)$ given in Figure~6. As before, $L$ will simultaneously denote 
the orthogonal block structure on $\Z_{36}$ consisting of coset-partitions of $Z_l, l \in L$.

In order to obtain $L$ as a poset block structure we start with the poset
$N=([4],\preceq)$ depicted in part (i) of Figure~13.
The dual lattice of ancestral subsets of $N$ has Hasse diagram shown in part (ii) of Figure~13.
This is indeed isomorphic to our lattice $L$.
Next, let us choose sets $X_1=[3]$, $X_2=[2]$, $X_3=[3]$ and $X_4=[2]$. We define the mapping
$f \colon X_1 \times X_2 \times X_3 \times X_4 \to \Z_{36}$ as
$$ (x_1,x_2,x_3,x_4) \mapsto 12x_1+18x_2+2x_3+9x_4 \pmod {36}.$$
The reader is invited to work out that $f$ is a bijection, and that $f$ is a week isomorphism from
the poset block structure defined by $N$ and the sets $X_i$ to our block structure $L$. \QED

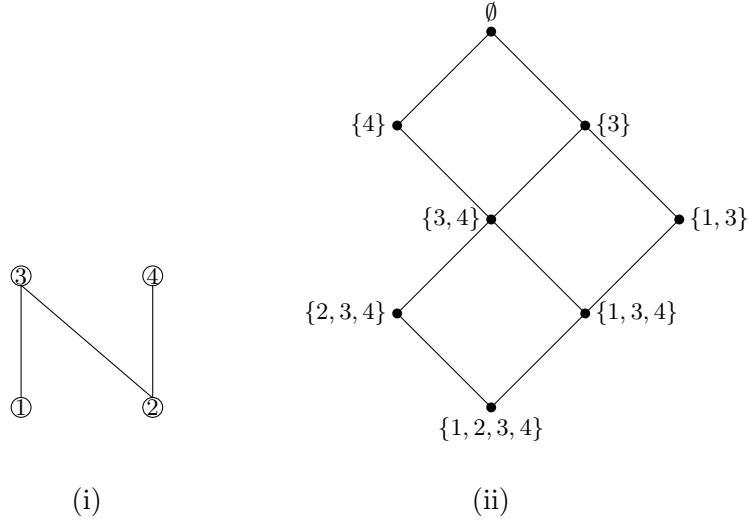
\begin{figure} 
\centering

\begin{tikzpicture}[scale=1.25] ---
\fill (0,0) circle (1.5pt);
\draw (0,0) node[below] {\footnotesize $\{1,2,3,4\}$};
\fill (-1,1) circle (1.5pt);
\draw (-1,1) node[left] {\footnotesize $\{2,3,4\}$};
\fill (1,1) circle (1.5pt);
\draw (1,1) node[right] {\footnotesize $\{1,3,4\}$};
\fill (0,2) circle (1.5pt);
\draw (0,2) node[left] {\footnotesize $\{3,4\}$};
\fill (2,2) circle (1.5pt);
\draw (2,2) node[right] {\footnotesize $\{1,3\}$};
\fill (1,3) circle (1.5pt);
\draw (1,3) node[right] {\footnotesize $\{3\}$};
\fill (-1,3) circle (1.5pt);
\draw (-1,3) node[left] {\footnotesize $\{4\}$};
\fill (0,4) circle (1.5pt);
\draw (0,4) node[above] {\footnotesize $\emptyset$};
\draw (0,0) -- (-1,1) -- (0,2) -- (1,1) -- (0,0);
\draw (1,1) -- (2,2) -- (1,3) -- (0,2);
\draw (1,3) -- (0,4) -- (-1,3) -- (0,2);
\draw (0,-1) node {\textrm (ii)};

\draw (-5,0) circle (3pt);
\draw (-5,0) node {\footnotesize $1$};
\draw (-3.6,0) circle (3pt);
\draw (-3.6,0) node {\footnotesize $2$};
\draw (-5,1.4) circle (3pt);
\draw (-5,1.4) node {\footnotesize $3$};
\draw (-3.6,1.4) circle (3pt);
\draw (-3.6,1.4) node {\footnotesize $4$};
\draw (-3.6,0.1) -- (-3.6,1.3);
\draw (-5,0.1) -- (-5,1.3);
\draw (-3.6,0.1) -- (-5,1.3);
\draw (-4.3,-1) node {\textrm (i)};
\end{tikzpicture}

\caption{Poset $N$ and the dual lattice of its ancestral subsets.}
\end{figure}

\end{exm}

\medskip

Now we are approaching the group-theoretical concept, crucial for the current presentation.
Let $\F$ be a poset block structure defined by a poset $(I,\preceq)$ and sets
$X_i$ $(i\in I)$. Recall that $A_i=\{j \in I \, \mid \, i \prec j\}$ is uncestral for all $i \in I$.
We set
$$ X^i = \prod_{j\in A_i}X_j = \prod_{i \prec j}X_j,$$
and $\pi^i$ for the {\em projection} of $X=\prod_{i\in I}X_i$ onto $X^i$.
The following construction can be found in \cite{BaiPRS83}.

\begin{defi}\label{def:GWP}
Let $(I,\preceq)$ be a poset, $X_i$ be a set ($i\in I$),
$|X_i| \ge 2$, and $K_i$ be a permutation group $K_i \le \sym(X_i)$.
The {\em generalized wreath product} $\prod_{(I,\preceq)}K_i$
defined by $(I,\preceq)$ and the groups $K_i$, is the complexus product
$$ 
P = \prod_{i \in I}P_i,
$$ 
where $P_i$ is the permutation representation of the group
$K_i^{X^i}$ on $X$ acting by the rule
$$ 
(x^f)_j =\left\{ \begin{array}{rr} x_j^{f(\pi^i(x))} & \hbox{if }i=j \\
                                      x_j &\hbox{if }i \ne j
                    \end{array}\right. , \; x=(x_j) \in X, \; f \in K_i^{X^i}.
$$
\end{defi}

\medskip

\noindent Clarification of the notation $f(\pi^i(x))$ follows below. 
We remark that, the fact that the above complexus product is indeed a group was proved in
\cite{BaiPRS83}. This construction has very interesting history, see Section~12.

Let $I=[r]=\{1,\ldots,r\}$ in Definition~\ref{def:GWP}. We write $x=(x_1,\ldots,x_r)$
for $x \in X=\prod_{i=1}^r X_i$. Every $f \in P$ is presented uniquely as the product
$f=f_1\cdots f_r$, where each $f_i \in P_i$.
Analogously to the ordinary wreath product (see 2.1), we shall also write
$f$ in the {\em table form}
$$ f=\big[ \; f_1(\pi^1(x)),\ldots,f_r(\pi^r(x)) \; \big].$$
By definition, $(x_1,\ldots,x_r)^f=\big(x_1^{f_1(\pi^1(x))},\ldots,x_r^{f_r(\pi^r(x))}\big).$
It is not hard to see that the group $P =\prod_{([r],\preceq)}K_i$ has order
\begin{equation}\label{eq4}
\big|\prod_{([r],\preceq)}K_i\big|=\prod_{i=1}^{r} |K_i|^{m_i},
\end{equation}
where $m_i=1$ if $\{i\} \in \anc((I,\preceq))$, and
$m_i=\prod_{j \in A_i}|X_j|$ otherwise.
The generalized wreath product gives back the ordinary direct and wreath product.
Namely, in case $r=2$ and the poset is an anti-chain the group $P=K_1 \times K_2$, and
if the poset is a chain with $1 \prec 2$, then  $P = K_2 \wr K_1$.

The following result about the automorphism group of a poset block
structure was proved by Bailey et al. (see \cite[Theorem A]{BaiPRS83}).
We say that a poset $(I,\preceq)$ satisfies the {\em maximal condition} if any subset
$J \subseteq I$ contains a maximal element.

\begin{thm}\label{thm:BPRS}
Let $(I,\preceq)$ be a poset having the maximal condition,
$X_i$ be a set of cardinality  at least $2$ for all $i \in I$, and
$\F$ be the poset block structure on $X$ defined by $(I,\preceq)$ and
the sets $X_i$. Then $\aut(\F)=\prod_{(I,\preceq)}\sym(X_i)$.
\end{thm}

Of course, if the set $I$ is finite, then $(I,\preceq)$ satisfies the maximal condition.
In particular, the above theorem applies to the orthogonal group block structures on $\Z_n$,
and hence we observe that their automorphism groups are certain generalized wreath products.
As an illustration of the above ideas, we determine once more the automorphism group of a rational
circulant graph, corresponding to our lattice $L$, in terms of generalized wreath product.

\begin{exm}\label{exm:graph}
Let $\Gamma$ be the rational circulant graph $\cay(\Z_{36},S)$, where

\begin{eqnarray*}
S &=& \{ 2,3,4,6,8,10,14,15,16,20,21,22,26,28,30,32,33,34 \} \\
  &=& (\Z_{36})_2 \cup (\Z_{36})_3 \cup (\Z_{36})_4 \cup (\Z_{36})_6.
\end{eqnarray*}

Because of Theorem~\ref{thm:W} the group $\aut(\Gamma) = \aut(\sgg{S}),$ where $\sgg{S}$ is the S-ring
over $\Z_{36}$ generated by $S$. S-ring $\sgg{S}$ is rational, hence by
Theorem~\ref{thm:M}, $\sgg{S} = \sg{ \un{Z_d} \mid d\in L}$ for a sublattice
$L$ of $L(36)$. After some simple reasonings we see that $L$ is our sublattice in Figure~5.
Thus 
$$
\aut(\Gamma)=\aut(\sgg{S})=\aut(L).
$$

As shown in Example~\ref{exm:n=36'''}, the orthogonal block structure $L$ is weekly isomorphic to
the poset block structure $\F$ defined by the poset $N=([4],\preceq)$ and sets
$X_i=[n_i], i \in \{1,\ldots,4\}$. Therefore, $\aut(L)$ is permutation isomorphic to the group
$\aut(\F)$. By Theorem~\ref{thm:BPRS}, the latter group $\aut(\F) = \prod_{N}S_{n_i}$
(we may get order once more using formula \eqref{eq4} as
$| \, \prod_{N}S_{n_i} \, |=(3!)^3 \cdot (2!)^6 \cdot 3! \cdot 2!=2^{11} \cdot 3^4.$)
\QED

\medskip

We converge with the consideration of our striking example.  Simultaneously, in principle, the main goals of the paper are fulfilled. Combination of all presented results implies that the automorphism groups of rational circulant graphs  are described by the groups as they appear in Theorem~\ref{thm:BPRS}.
Nevertheless, at this stage we are willing to justify much more precise formulation,
as it is presented in the main Theorem~\ref{thm:main}, as well as to provide its self-contained proof.
\end{exm}

\section{Proof of Theorem~\ref{thm:main}}

Let $P=([r],\preceq)$ be a poset, and $n_1,\ldots,n_r$ be in $\N$ such that
$n_i \ge 2$ for all $i \in \{1,\ldots,r\}$. We denote by
$\pbs(P;n_1,\ldots,n_r)$ the poset block structure defined by $P$ and the sets
$[n_i]$. We recall that $P=([r],\preceq)$ is increasing if
$i \preceq j$ implies $i \le j$ for all $i,j \in [r]$.

\medskip

The final step toward Theorem~\ref{thm:main} is the following statement.

\begin{prop}\label{prop:step3} \
\begin{enumerate}[(i)]
\item Let $P=([r],\preceq)$ be an increasing poset and
$n_1,\ldots,n_r$ be in $\N$ satisfying
\begin{enumerate}[(a)]
\item $n=n_1 \cdots n_r$,
\item $n_i \ge 2$ for all $i \in \{1,\ldots,r\}$,
\item $(n_i,n_j)=1$ for all $i,j \in \{1,\ldots,r\}$ with $i \not\preceq j$.
\end{enumerate}
Then  $\pbs(P;n_1,\dots,n_r)$ is weakly isomorphic to an
orthogonal group block structure on $\Z_n$.
\item Let $\F$ be an orthogonal group block structure on $\Z_n$.
Then exists an increasing poset $P=([r],\preceq)$ and
$n_1,\ldots,n_r$ in $\N$ satisfying (a)-(c) in (i) such that $\F$ is
weakly isomorphic to $\pbs(P;n_1,\ldots,n_r)$.
\end{enumerate}
\end{prop}

To settle the proposition we first prove two preparatory lemmas.
For $J \subset [r]$ we set the notation $\ov{J}=[r] \setminus J$.

\begin{lem}\label{lem:help1}
Let $P=([r],\preceq)$ be an increasing poset and
$n_1,\ldots,n_r$ be in $\N$ satisfying (a)-(c) in (i) of Proposition~\ref{prop:step3}.
Then the set $L=\big\{ \prod_{j\in \ov{J}}n_j \, \mid \,  J \in \anc(P)\big\}$ is a sublattice of $L(n)$.
\footnote{If $J \in \anc(P)$ is the whole set $[r]$, then we set $\prod_{j\in \ov{J}}n_j=1$.}
\end{lem}

\proof We prove the lemma by induction on $r$. If $r=1$ then $L=\{1,n\}$.
Suppose that $r>1$ and let $n'=n_1 \cdots n_{r-1}$.
Let $P'=([r-1],\preceq)$ be the poset on $[r-1]$ induced by $\preceq$.
The induction hypothesis applies to $P'$ and numbers $n_1,\ldots,n_{r-1}$.
Thus we get sublattice $L_2$ of $L(n')$ as
$$ L_2 = \big\{ \, \prod_{j\in \ov{J}}n_j \, \mid \, J \in \anc(P') \, \big\}. $$
Here by $\ov{J}$ we mean the complement of $J$ in $[r-1]$.

Since $P$ is increasing, $r$ is a maximal element in $P$.
Thus for any $J \subseteq [r-1]$,
\begin{equation}\label{eq7}
J \in \anc(P') \iff J \cup \{r\} \in \anc(P).
\end{equation}
Let $J_*=\{j \in [r-1] \, \mid \, j \not\preceq r \}$. Then $J_* \in \anc(P)$. Further, for any
$J \subseteq [r-1]$,
\begin{equation}\label{eq8}
J \in \anc(P) \iff J \in \anc(P') \textrm{ and } J \subseteq J_*.
\end{equation}

Put $d=\prod_{j \in \ov{J_*}}n_j$. Clearly, $d \in L_2$.
Then $n'/d=\prod_{j \in [r-1],j \not\preceq r}n_j$, hence condition (c)
in (i) of Proposition~\ref{prop:step3} implies that $n'/d \wedge n_r=1$.
Thus we can use Definition~\ref{def:crestedp2} to form the crested product $L_1 \otimes_d L_2,$ where
$L_1=\{1,n_r\}$. Then
\begin{eqnarray*}
L_1 \otimes_d L_2 &=&
\big\{ \,  l_1l_2 \, \mid \, l_1=1, l_2\in L_2, \text{ or } l_1 \in L_1,
l_2 \in L_2 \text{ with } d \mid l_2 \, \big\} \\
&=& L_2 \cup \big\{ \, n_r l_2 \, \mid \, l_2 \in L_2 \text{ with } d \mid L_2 \, \big\}.
\end{eqnarray*}

Now, we use \eqref{eq7} and \eqref{eq8} to find
\begin{eqnarray*}
L &=& \Big\{\prod_{j\in \ov{J}}n_j \, \mid \, J \in \anc(P) \; \textrm{ and }
r \in J \Big\} \cup \Big\{ \prod_{j\in \ov{J}}n_j \, \mid \, J \in
\anc(P) \; \textrm{ and } r \notin J \Big\} \\
  &=& L_2 \cup \big\{ \, l_2n_r \, \mid \, l_2 \in L_2 \text{ with } d \mid l_2 \, \big\}=L_1 \otimes_d L_2.
\end{eqnarray*}
Thus $L$ is a sublattice of $L(n),$ as required. \QED

We show next the converse to Lemma~\ref{lem:help1}.

\begin{lem}\label{lem:help2}
Let $L$ be a sublattice of $L(n), n \ge 2$ such that $1,n \in L$.
Then $L = \big\{ \prod_{j\in \ov{J}}n_j \, \mid \,  J \in \anc(P)\big\},$ where
$P=([r],\preceq)$ is an increasing poset and
$n_1,\ldots,n_r$ are in $\N$ satisfying (a)-(c) in (i) of Proposition~\ref{prop:step3}.
\end{lem}

\proof We proceed by induction on $n$. The statement is clear if $L=\{1,n\}$.
Suppose $L \ne \{1,n\}$, and let $m$ be a maximal element in the poset
induced by $L \setminus \{n\}$. Induction applies to sublattice $L\subl{m}$, and we can write
$$ L\subl{m} = \big\{ \, \prod_{j \in \ov{J}}n_j \, \mid \, J \in \anc(P') \, \big\}$$
with a suitable poset $P'=([r-1],\preceq)$ and numbers $n_1,\ldots,n_{r-1}$
in $\N$.
Now, let $s$ be the smallest number in the set $L \setminus L\subl{m}$.
Since $m \wedge s \in L\subl{m},$ we have a subet $J_* \in \anc(P')$ for which
$m \wedge s=\prod_{j \in \ov{J}_*}n_j.$
Define the poset $P$ on $[r]$ as the extension of $P'$ to $[r]$ by setting
$r \not\preceq x$ for all $x \in [r-1]$, and $$ x \preceq r \iff x \notin J_*.$$
We claim that $P$ is the required poset and $n_1,\ldots,n_{r-1},n_r=n/m$ are
the required numbers.

First, poset $P$ is obviously increasing, $n_1 \cdots n_r=n$, and
$n_i \ge 2$ for all $i \in [r]$. Let $i,j \in [r]$ with $i \not\preceq j$. It is clear that
$n_i \wedge n_j=1$ if $j \ne r$. Let $j=r$. Then
$$ n_r = \frac{n}{m} = \frac{s}{m \wedge s}, \text{ and }
\prod_{j \not\preceq r}n_j = \prod_{j \in J_*}n_j = \frac{m}{m \wedge s}.$$
This shows that $n_i \wedge n_r=1$ holds as well, and so $n_1,\ldots,n_r$ satisfy
(a)-(c) in (i) of Proposition~\ref{prop:step3}.

\medskip

By \eqref{eq7} and \eqref{eq8},
\begin{eqnarray*}
\Big\{ \, \prod_{j \in \ov{J}}n_j \, \mid \, J \in \anc(P) \, \Big\} &=&
\Big\{ \, \prod_{j\in \ov{J}}n_j \, \mid \, J \in \anc(P) \, \text{ and }
r \in J \, \Big\} \cup \\
& & \Big\{ \, \prod_{j\in \ov{J}}n_j \, \mid \, J \in \anc(P) \, \text{ and }
r \notin J \, \Big\} \\
&=& L\subl{m} \cup \big\{ \, xn_r \, \mid \, x \in (L\subl{m})\supl{m \wedge s} \, \big\}.
\end{eqnarray*}
Now, we use Lemma~\ref{lem:M} to conclude
$$ L\subl{m} \cup \Big\{ \, x \, \frac{s}{m \wedge s} \, \mid \, x \in (L\subl{m})\supl{m \wedge s} \, \Big\}
=L\subl{m} \cup (L \setminus L\subl{m})=L. $$ \QED

\noindent{\sc Proof of Proposition~\ref{prop:step3}.} \quad
Let $P=([r],\preceq)$ be an increasing poset and
$n_1,\ldots,n_r$ be in $\N$ satisfying (a)-(c) in (i) of Proposition~\ref{prop:step3}.
Let $L = \big\{ \prod_{j\in \ov{J}}n_j \, \mid \,  J \in \anc(P)\big\}$ be the
sublattice of $L(n)$.
In view of Lemmas~\ref{lem:help1} and ~\ref{lem:help2} it remains to prove that
$\pbs(P;n_1,\ldots,n_r)$ is weakly isomorphic to the orthogonal group block structure on $\Z_n$ defined by $L$.

\medskip

Let $J \in \anc(P)$, $J \ne [n]$, and $x_j,y_j \in [n_j]$ for each $j \in \ov{J}$.
We claim that
\begin{equation}\label{eq9}
 \sum_{j \in \ov{J}}\Big( \prod_{i \in [r],i \not\preceq j}n_i \Big) x_j \equiv
   \sum_{j \in \ov{J}}\Big( \prod_{i \in [r],i \not\preceq j}n_i \Big) y_i \pmod
{n} \; \implies \; x_j=y_j \textrm{ for each }j \in \ov{J}.
\end{equation}
We proceed by induction on $r$. Let $r=1$. Then $J=\emptyset$, the assumption in \eqref{eq9}
reduces to $x_1 \equiv y_1 \pmod n$, and from this $x_1=y_1$.

Let $r>1$. Let $n' = n/n_r$ and $P'$ be
the poset induced by $[r-1]$. Let $J \in \anc(P)$. First, let $r \in J$, and put
$J'=J \setminus \{r\}$. By \eqref{eq7}, $J' \in \anc(P')$.
The assumption in \eqref{eq9} can be rewritten in the form
$$\sum_{j \in \ov{J'}}\Big( \prod_{i \in [r-1],i \not\preceq j}n_i \Big)n_rx_j \equiv
\sum_{j \in \ov{J'}}\Big( \prod_{i \in [r-1],i \not\preceq j}n_i \Big)n_ry_i
\pmod {n}.$$
From this
$$\sum_{j \in \ov{J'}}\Big( \prod_{i \in [r-1],i \not\preceq j}n_i \Big) x_j \equiv
\sum_{j \in \ov{J'}}\Big( \prod_{i \in [r-1],i \not\preceq j}n_i \Big)y_i \pmod {n'}, $$
and hence, by induction, $x_j=y_j$ for each $j \in \ov{J'}=\ov{J}$, and \eqref{eq9} holds.
Second, let $r \notin J$. Put $n_r^*=\prod_{i\in [r],j \not\preceq r}n_j$.
Notice that $n_r \wedge n_r^*=1$ (see (c) in (i) of Proposition~\ref{prop:step3}).
The assumption in \eqref{eq9} can be rewritten as
$$\sum_{j \in \ov{J},j \ne r}
\Big( n_r \prod_{i \in [r-1],i \not\preceq j}n_i \Big) x_j + n_r^*x_r \equiv
\sum_{j \in \ov{J}, j \ne r}
\Big( n_r \prod_{i \in [r-1], i \not\preceq j}n_i \Big) y_i + n_r^*y_r \pmod {n}. $$
From this $n_r^*(x_r - y_r) \equiv 0 \pmod {n_r}$. And as $n_r \wedge n_r^*=1$, $x_r=y_r$.
By \eqref{eq8}, $J \in \anc(P')$. Regarded $J$ as an ancestral subset of $P'$, we find
$$\sum_{j \in \ov{J}}\Big( \prod_{i \in [r-1], i \not\preceq j}n_i \Big) x_j \equiv
\sum_{j \in \ov{J}}\Big( \prod_{i \in [r-1],i \not\preceq j}n_i \Big) \pmod {n'}. $$
Thus, by induction, $x_j=y_j$ for each $j \in [r-1] \setminus J$, and so \eqref{eq9} holds.

\medskip

Let $X=[n_1] \times \cdots \times [n_r]$.
Define the mapping
$$f \colon X \to \Z_n, \; (x_i) \mapsto \sum_{i=1}^{r}\big(\prod_{j \in [r],
j \not\preceq i}n_j \big) x_i \pmod n.$$
We claim that $f$ is a weak isomorphism from $\pbs(P;n_1,\ldots,n_s)$ to $L$.
First, that $f$ is a bijection can be seen from \eqref{eq9} by substituting $J=\emptyset$.
Let $J \in \anc(P)$, and fix an element $(x_i)=(x_1,\ldots,x_r) \in X$.
Recall that the class of $\Pi(J)$ containing $(x_i)$ is the set
$$\big\{ \; (y_i) \in X \; \mid \; x_j=y_j \textrm{ \ for all \ }i \in J \; \big\}.$$
Put $m=\sum_{j \in J}\big(\prod_{i \in [r],
i \not\preceq j}n_i \big) x_j$ in $\Z_n$.
Then $f$ maps the above class to the set
$$ m+\Big\{ \; \sum_{j \in \ov{J}}\big(\prod_{i \in [r],
i \not\preceq j}n_i \big) y_j \; \mid \; j \in \ov{J}, \; y_j \in [n_j] \; \Big\}.$$
Observe that $i \not\preceq j$ for any $j \in \ov{J}$ and $i \in J$.
Thus the product $\prod_{j \in J}n_j$ divides the numbers in the above set, and hence
$$m+\Big\{ \; \sum_{j \in \ov{J}}\big(\prod_{i \in [r],
i \not\preceq j}n_i \big) y_j \; \mid \;  j \in \ov{J}, \; y_j \in [n_j] \; \Big\}
\subseteq m + \big\langle \; \prod_{j \in J}n_j \; \big\rangle=m + Z_d, \; d=
{\prod_{j\in \ov{J}}n_j}.$$
Since $f$ is a bijection, equality follows, and so $f$ maps the partition $\Pi(J)$ to the
coset-partition $F_{Z_d}$. This completes the proof of the proposition. \QED



\bigskip

\noindent{\sc Proof of Theorem \ref{thm:main}.}

\medskip

$(i) \Rightarrow (ii)$ \
Let $\cay(\Z_n,S)$ be a rational circulant graph with $G=\aut(\cay$ $(\Z_n,S))$.
By Corollary~\ref{cor:step2}, $G=\aut(\F)$, where $\F$ is an orthogonal group block structure on $\Z_n$.
By (ii) of Proposition~\ref{prop:step3}, $\F$ is weakly isomorphic to the poset block structure 
$\pbs(P;n_1,\ldots,n_r)$ for suitable poset $P=([r],\preceq)$ and numbers $n_1,\ldots,n_r$.
Theorem~\ref{thm:BPRS} gives that $G$ is permutation isomorphic to $\Pi_P S_{n_i}$.

\medskip

$(ii) \Rightarrow (i)$ \ Let $G=\Pi_P S_{n_i}$, where
$P=([r],\preceq)$ is an increasing poset and $n_1,\ldots,n_r$ are in $\N$ satisfying
(a)-(c) in (ii) of Proposition~\ref{prop:step3}.
Because of Theorem~\ref{thm:BPRS} the group $G$ equals the automorphism group of the
poset block structure $\pbs(P;n_1,\ldots,n_r)$.
By (i) of Proposition~\ref{prop:step3}, $\pbs(P;n_1,\ldots,n_r)$ is weakly isomorphic to
an orthogonal group block structure $\F$ on $\Z_n$, hence $G$ is permutation isomorphic to $\aut(\F)$.
Finally, Corollary~\ref{cor:step2} shows that there is a
rational circulant graph $\cay(\Z_n,S)$ such that $\aut(\cay(\Z_n,S))=\aut(\F)$. \QED

\section{Miscellany}

We conclude the paper by a collection of miscellaneous topics related to rational
circulant graphs and their automorphisms.

\subsection{Enumeration of rational circulant graphs.} 

Let $\cay(\Z_n,S)$ be a rational graph. By Theorem~\ref{thm:BM}, $S$ follows to be
the union of some of the sets
$$ (\Z_n)_d=\{ x \in \Z_n \, \mid \, \gcd(x,n)=d \},$$
where $d$ is a divisor of $n$. Conversely, any set in such a form is a
connection set of a rational graph. In particular, up to isomorphism, we have at most
$2^{\tau(n)-1}$ rational Cayley graphs (without loops) over $\Z_n$.

To investigate, which of these graphs are pairwise non-isomorphic, we refer to the following Zibin's conjecture
for arbitrary circulant graphs, which follows easily from the results in \cite{Muz94}
(see also \cite{MuzP99} and \cite[Theorem~5.1]{MuzKP01}).

\begin{thm}\label{thm:MP} {\bf (Zibin's conjecture)}
Let $\cay(\Z_n,S)$ and $\cay(\Z_n,R)$ be two
isomorphic circulant graphs. Then for each $d \mid n$ there exists a multiplier $m_d \in \Z_n^*$
such that $S_d^{(m_d)}=R_d$.
\end{thm}

Here for arbitrary subset $S \subseteq \Z_n,$ we define $S_d=S \cap (\Z_n)_d$.

\begin{cor}
Let $\cay(\Z_n,S)$ and $\cay(\Z_n,R)$ be two rational circulant graphs.
Then these are isomorphic if and only if $S_d=R_d$ for all $d \mid n$.
Moreover, for each $d \mid n$ the common set $S_d=R_d$ is equal to $\emptyset$ or to 
$(\Z_n)_d$. 
\end{cor}

\begin{cor}
The number of non-isomorphic rational circulant graphs (without loops) of order $n$
is $2^{\tau(n)-1}$, where $\tau(n)$ is the number of positive divisors of $n$.
\end{cor}

We should mention that the above statement was given as a conjecture by
So \cite[Conjecture~7.3]{So05}.

\medskip

\noindent {\bf Remark 1.} \ We refer to the sequence A100577 (starting with $1,2,2,4,2,8,2,8,4,8,$ $2,32$) 
in the famous Sloane's On-Line Encyclopedia of Integer Sequences, see \cite{OEIS}, which consists of the 
numbers $2^{\tau(n)-1}, n \in \N$.

\medskip

\noindent {\bf Remark 2.} \ For a given $n$ let $X$ be an arbitrary subset of the set $L(n) \setminus \{n\}$.
Let $S=\cup_{d\in X} (\Z_n)_d,$ $\Gamma=\cay(\Z_n,S)$. Clearly, $\Gamma$ is a presentation of an arbitrary rational
circulant with $n$ vertices and the group $\aut(\Gamma)$ is completely defined by a suitable sublattice $L$ of 
the lattice $L(n)$. In Example 9.5 for $X$ presented there the corresponding lattice $L$ coincides with our 
striking sublattice.  

A question of elaboration of a simple procedure to recognize $L$ from an arbitrary subset $X$ is of a definite 
independent interest, though it is out of the scope in the current text.

\subsection{Association schemes.} 

Though we have managed to arrange the main line of the presentation
without the evident use of association schemes, it is now time to consider explicitly this concept.

Let $X$ be a nonempty finite set, and let
$\Delta_X$ denote the {\em diagonal relation} on $X$, i. e.,
$\Delta_X=\{(x,x) \mid x \in X\}$.
For a relation $R \subseteq X \times X$, its {\em transposed} $R^t$ is defined by
$R^t=\{(y,x) \mid (x,y) \in R\}$. For a set $\{R_0,R_1,\ldots,R_d\}$
of relations on $X$ the pair $\X=(X,\{R_0,\ldots,R_d\})$ is called an
{\em association scheme} on $X$ if the following axioms hold (see \cite{BanI84}):
\begin{enumerate}[({AS}1)]
\item $R_0=\Delta_X$, and $R_0,R_1,\ldots,R_d$ form a partition of $X \times X$.
\item For every $i \in \{0,\ldots,d\}$ there exists $j \in \{0,\ldots,d\}$ such that
$R_i^t=R_j$.
\item For every triple $i,j,k \in  \{0,\ldots,d\}$ and for $(x,y) \in R_k$, the number,
denoted by $p^{k}_{i,j}$, of elements $z \in X$ such that $(x,z) \in R_i$ and
$(z,y) \in R_j$ does not depend
on the choice of the pair $(x,y) \in R_k$.
\end{enumerate}
The relations $R_i$ are called the {\em basic relations} of $\X$, the
corresponding graphs $(X,R_i)$ the {\em basic graphs} of $\X$.
The automorphism group of $\X$ is the permutation group
$$\aut(\X):=\bigcap_{i=0}^{r}\aut((X,R_i)).$$

\medskip

Let $\F$ be an orthogonal block structure on $X$.
For $F \in \F$,  define the {\em color relation} $C_F$ on $X$ as
$$(x,y) \in C_F \iff F=\bigwedge \big\{ E \in \F \, \mid (x,y) \in R_E \big\}.$$

It is immediately clear that for each $F \in \F$, $C_F$ is a symmetric relation and the set of relations $\{C_F, F \in \F \}$ forms a partition of $X \times X$.  It turns out that we can claim more.

The relational system $(X,\{C_F \, \mid \, F \in \F\})$ is a symmetric association
scheme on $X$ (see \cite[Theorem 4]{Bai96}). Recall that in a {\em symmetric} association scheme
$R_{i}^{t} = R_i$ for all $i \in \{0,1, \ldots d\}$. This we denote by $\co(\F)$.
Observe that, if $\F$ is the orthogonal group block structure on $\Z_n$ given in
Proposition~\ref{prop:autA=autF}, then the color relations are $\widehat{R_l}$
defined in the proof of Proposition~\ref{prop:autA=autF}.

\subsection{Schurity of rational S-rings over cyclic groups.} 

Let $\A$ be a rational S-ring over $\Z_n$, and $\F$ be the corresponding
orthogonal group block structure on $\Z_n$. Recall that $\A$ is Schurian if
$\A=V(\Z_n,\aut(\A)_e)$. It is not hard to see that this is equivalent to saying that
the basic relations of the association scheme $\co(\F)$ are the $2$-orbits
of $\aut(\A)$, the latter group is the same as $\aut(\F)=\aut(\co(\F))$).

The following result is due to Bailey et al. (see \cite[Theorem C]{BaiPRS83}),
which in particular also answers the Schurity of
rational S-rings over $\Z_n$ in the positive.

\begin{thm}
Let $(I,\preceq)$ be a finite poset, $X_i$ be a finite set of cardinality at
least $2$ for each $i \in I$, $X=\prod_{i\in I}X_i$, and $\F$ be the poset block
structure on $X$ defined by $(I,\preceq)$ and the sets $X_i$. Then the association
scheme $\co(\F)$ is Schurian, i. e., $\co(\F)=(X,\orb(\aut(\co(\F))))$.
\end{thm}

\begin{cor}
Every rational S-ring over $\Z_n$ is Schurian.
\end{cor}

Let us remark that it has been conjectured that all S-rings over the
cyclic groups $\Z_n$ are Schurian (known also as the Schur-Klin conjecture).
The conjecture was denied recently by Evdokimov and Ponomarenko \cite{EvdP02}.

\subsection{Simple reduction rules.} 

We say that {\em simple reduction rules} apply to the group $\Z_n$ if every 
sublattice $L$ of $L(n)$ such that $1,n \in L,$ is obtained from trivial lattices 
via an iterative use of reduction rules 1 and 2.
For instance, simple reduction rules apply to $\Z_{12}$,
but not to $\Z_{36}$ (see the striking example).
We already discussed informally which are the orders $n$ that simple reduction rules apply to $\Z_n$.

\medskip

The question is closely related to simple block structures
introduced by Nelder \cite{Nel65}. Next we recall shortly the definition and some
properties following \cite{Bai96}.

For $i=1,2$, let $F_i$ be a partition of a set $X_i$. Define the partition
$(F_1,F_2)$ of $X_1 \times X_2$ by setting the corresponding equivalence relation
$R_{(F_1,F_2)}$ as
$$
\big( \; (x_1,x_2),(y_1,y_2) \; \big) \in R_{(F_1,F_2)} \iff (x_1,y_1) \in R_{F_1} \land
(x_2,y_2) \in R_{F_2}.
$$
Let $\F_i$ be a block structure on $X_i$ ($i=1,2$). Their {\em crossing product} is
the block structure on $X_1 \times X_2$ defined by
$$
\F_1 * \F_2 = \big\{ (F_1,F_2) \, \mid \, F_1 \in \F_1, \; F_2 \in \F_2 \big\},
$$
and their {\em nesting product} is the block structure on $X_1 \times X_2$ defined by
$$
\F_1 / \F_2 = \big\{ (F_1,U_2) \, \mid \, F_1 \in \F_1 \big\} \;
\cup \; \big\{ (E_1,F_2) \, \mid \, F_2 \in \F_2 \big\}.
$$
For the automorphism groups we have
$\aut(\F_1 * \F_2)=\aut(\F_1) \times \aut(\F_2)$, and
$\aut(\F_1 / \F_2)$ $=\aut(\F_1) \wr \aut(\F_2)$.
Note that, if in addition both $\F_i$ are orthogonal, then so are
$\F_1 * \F_2$ and $\F_1 / \F_2$. The {\em trivial block} structure on a set $X$ is
the one formed by the partitions $E_X$ and $U_X$. The {\em simple} orthogonal block structures
are defined recursively as follows:
\begin{itemize}
\item Every trivial block structure is simple of depth $1$.
\item If for $i=1,2$, $\F_i$ is a simple orthogonal block structures of depth $s_i$ on a
set $X_i$, $|X_i| \ge 2$, then $\F_1 * \F_2$ and $\F_1 / \F_2$ are simple orthogonal block
structures of depth $s_1+s_2$.
\end{itemize}
Clearly, if $\F$ is simple, then $\aut(\F)$ is obtained
using iteratively direct or wreath product of symmetric groups.
The equivalence follows.

\begin{cor}
Simple reduction rules apply to $\Z_n$ if and only if every orthogonal group block structure on
$\Z_n$ is simple.
\end{cor}

\begin{cor}
Simple reduction rules apply to $\Z_n$ if and only if $n=pqr$, or $n=p^eq$, or $n=p^e,$
where $p,q$ and $r$ are distinct primes.
\end{cor}

\proof In view of the previous corollary we only need to check if there exists an orthogonal
group block structure $\F$ on $\Z_n$ which is not simple. By Proposition~\ref{prop:step3},
$\F$ is weakly isomorphic to $\pbs(P;n_1,\ldots,n_r),$ where $P=([r],\preceq)$ is a
non-increasing poset with suitable weights $n_i$.
Let $N$ be the poset given in part (i) of Figure~13.
It is proved that $\F$ is not simple if and only if $P$ contains a subposet isomorphic to $N$
(see \cite[pp. 64]{Bai96}). Let $m_i$, $1\le i \le 4$, be the weights of this subposet.
Then $m_1m_2m_3m_4 \ | \ n$, and hence $n \ne pqr$ for distinct primes $p$, $q$ and $r$.
Let $n=p^eq$ or $n=p^e$. Then $q$ appears as a factor in at most one of the numbers $m_i$,
and so $(m_1,m_2)>1$ or $(m_3,m_4)>1$. This contradicts (c) in (ii) of
Theorem~\ref{thm:main}. These yield implication `$\Leftarrow$' in the statement.

For implication `$\Rightarrow$' assume that $n$ is none of the numbers $pqr, p^eq$, or $p^e$, 
where $p,q$ and $r$ are distinct primes. 
We leave for the reader to check that in this case it is possible
to assign weights $n_i$ to $N$ satisfying (a)-(c) in (ii) of Theorem~\ref{thm:main}.
The arising orthogonal group block structure on $\Z_n$ is therefore not simple, and
by this the proof is completed. \QED

\section{Historical digest}

This paper objectively carries certain interdisciplinary features.
Indeed, main concepts discussed by us may be attributed to such areas as association schemes,
S-rings, group theory, design of statistical experiments,
spectral graph theory, lattice theory, etc.
While for the authors there exists an evident natural impact of ideas borrowed from
many diverse areas, it is difficult to expect similar experience from each interested reader.
Nevertheless, at least brief acquaintance with the roots of the many facets of
rational circulants, may create an extra helpful context for the reader. This is
why we provide in the final section a digest of
historical comments. 
We did not try to make it comprehensive, hoping to come
once more in a forthcoming paper to discuss the plethora  of all detected lines with more detail.

\subsection{Schur rings.}

The concept of an S-ring goes back to the seminal paper of Schur \cite{Sch33},
the abbreviation S-ring was coined and used by R.~Kochend\"orfer and H.~Wielandt.
For a few decades S-rings were used exclusively in permutation group
theory in framework of very restricted area of interests.  Books \cite{Sco64,DixM96}
provide a nice framework, showing evolution in attention of modern experts to
this concept. (Indeed, while S-rings occupy a significant
position in \cite{Sco64}, the authors of \cite{DixM96} avoid to use the term itself,
though still present the background of the classical applications of S-rings to so-called
B-groups, B stands for Burnside.)

On the dawn of algebraic graph theory, the interest to S-rings was revived due
to their links with graphs and association schemes, admitting a regular group as a
subgroup of the full automorphism group. In this context paper \cite{Cha67} by
C.~Y.~Chao definitely deserves credit for pioneering contribution.
More evident combinatorial applications of S-rings stem from \cite{Pos74,KliP81}.
Tendencies of modern trends for attention to the use of S-rings in graph theory
still are not clear enough. On one hand, a number of experts do not even try to
hide their fearful feelings toward S-rings, regarding their
ability to avoid ``heavy use of Schur rings'' (see \cite{Jos95}) as a
definite positive feature of their presentation.
On other hand, S-rings form a solid part of a background for high level
monographs, though under alternative names like translation association
scheme \cite{BroCN89} or blueprint \cite{Bai04}.

\subsection{Schur-rings over $\Z_n$.}

Practical application of established theory by Schur \cite{Sch33} originally was 
consideration of S-rings over finite cyclic groups.
As a consequence, he proved that every primitive overgroup
of a regular cyclic group of composite order $n$ in symmetric group $S_n$ is doubly transitive.
Further generalizations of this result are discussed in \cite{Wie64}.
Nowadays, the group theoretical results of such flavor are obtained with the aid of
classification of finite simple groups (CFSG), see e. g. \cite{Li03}.
Schur himself did not try to describe all S-rings over $\Z_n$.
First such serious attempt was done by P\"oschel \cite{Pos74} on suggestion of
L.~A.~Kalu\v{z}nin, disciple of Schur. In \cite{Pos74} all S-rings over cyclic groups of
odd prime-power order were classified. Classification of S-ring over group $\Z_{2^e}$
was fulfilled by joint efforts of Ja.~Ju.~Gol'fand, M.~H.~Klin, N.~L.~Naimark and R.~P\"oschel (1981-1985),
see references in \cite{MuzKP01,Kov05}.
First attempts of description of automorphism
groups of circulants of order $n$, their normalizers in $S_n$ and, as a result, a solution
of isomorphism problem for circulants can be traced to \cite{KliP81}.
K.~H.~Leung, S.~L.~Ma and S.~H.~Man reached complete recursive description of S-rings
over $\Z_n$ in \cite{LeuM90,LeuM96,LeuM98}.
An alternative approach was established by Muzychuk, see e.g. \cite{Muz94,Muz95}.
The results of Leung and Ma were rediscovered by S.~A.~Evdokimov and I.~N.~Ponomarenko
\cite{EvdP02}. In fact, in \cite{EvdP02} a much more advanced result was presented:
evident description of infinite classes of non-Schurian S-rings over $\Z_n$.

In 1967 A.~\'{A}d\'{a}m \cite{Ada67} posed a conjecture:
two circulants of order $n$ are isomorphic if and only if they are conjugate with the
aid of a suitable multiplier from $\Z_n^*$. A number of mathematicians more or less
immediately presented diverse counterexamples to this conjecture.
Nevertheless, a more refined question was formulated: for which values of $n$ the conjecture
is true, see \cite{Pal87} and references in it.
A complete solution of this problem was given in \cite{Muz97}.
Later on Muzychuk \cite{Muz04} provided a necessary and sufficient condition for
two circulants of order $n$ to be isomorphic.  This monumental result
(as well as previous publications) of Muzychuk is based on skillful
combination of diverse tools, including deep use of S-rings.
Schur rings were also used for the analytical enumeration of circulant graphs,
see \cite{KliLP96,LisP00}.
Current ongoing efforts for the description of the automorphism groups of circulant
graphs are also based on the use of S-rings. For $n$ equal to odd prime-power and
$n=2^e$ the problem is completely solved, see \cite{Kli93,Kli94,Kov02,KliP,Kov}.
A polynomial time algorithm which returns the automorphism group of an arbitrary
circulant graph was recently constructed by Ponomarenko \cite{Pon06}.

For about four decades investigation of Schur rings over cyclic groups is serving for
generation of mathematicians as a challenging training polygon in development of algebraic
graph theory. This supports the author's enthusiasm to further promote combinatorial
applications of S-rings and to expose this theory to a wider audience.

\subsection{Rational S-rings and integral graphs.} 

Original name coined by Schur was S-ring of traces. It seems that
Wielandt \cite{Wie64} was the first who suggested to use adjective rational
for this class of S-rings.  The complete rational S-ring $\A_n$ over
$\Z_n$ appears as the transitivity module of the holomorph of $\Z_n$, which is isomorphic to
$\Z_n \rtimes \Z_n^*$.  Its basis quantities are orbits of the multiplicative action of
$\Z_n^*$ on $\Z_n$.
It was already Schur who noticed that in a similar way $\Z_n^*$ acts multiplicatively
on an arbitrary finite abelian group $H$ of exponent $n$.  Thus also in this case it
is possible to consider the transitivity module of
$H \rtimes \Z_n^*$.  The resulting S-ring is exactly the complete rational S-ring over $H$.
W.~G.~Bridges and R.~A.~Mena rediscovered in \cite{BriM79} (in a different context) the algebra $\A_n$
and exposed a lot of its significant properties.
Only later on, in \cite{BriM82}, they realized (due to hint of E.~Bannai)
existence of links of their generalization of $\A_n$ for arbitrary finite abelian groups
with the theory of S-rings. A crucial contribution, exploited in \cite{BriM79,BriM82},
was the use of the group basis in the complete rational S-ring over $H$.
Implicitly or explicitly the algebras $\A_n$ and $V(H,\Z_n^*)$ were investigated later on
again and again, basing on diverse motivation see e.g. \cite{RaoRS84,Fert85,Gol85,BanS02}.

As was mentioned, Muzychuk's classification of rational S-rings over $\Z_n$ \cite{Muz93}
forms a cornerstone for the background of the current paper. In turn, solutions for
two particular cases, that is $n$ is a prime-power \cite{Pos74} and $n$ is square-free
\cite{Gol85} created a helpful starting context for Muzychuk.
Essential tools exploited in \cite{Muz93} are use of group basis and possibility to
work with so-called pseudo-S-rings (those which do not obligatory include $\un{e}$ and
$\un{H}$). In fact, pseudo-S-rings were used a long time ago by Wielandt \cite{Wie64}.
This, in conjunction with the classical techniques of Schur ring theory, allows
to obtain transparent proofs of main results.
Zibin's conjecture (and its particular case Toida's conjecture) were also proved
in \cite{MuzKP01} with the aid of S-rings based on further results of Muzychuk.
An alternative approach developed in \cite{DobM02} depends on the use of CFSG.

\medskip

F.~Harary and A.~J.~Schwenk \cite{HarS74} suggested to call a graph $\Gamma$ 
{\em integral} if every eigenvalue of $\Gamma$ is integer. Since their pioneering paper a lot of
interesting results about such graphs were obtained. A very valuable survey appears in
\cite[Chapter 5]{PetR01}.  More fresh results are discussed in \cite{Wan05}.
It was proved in \cite{AhmABS09} that integral graphs are quite rare, that is, only a
fraction of $2^{-\Omega(n)}$ of the graphs on $n$ vertices have an integral spectrum.
Recent serious applications of integral graphs for designing the network topology of
perfect state transfer networks (see e.g. references in \cite{AhmABS09}) imply new wave
of interest to these graphs. In the context of the current paper, our interest
to integral graphs is strictly restricted by regular graphs.  A significant source of
regular integral graphs is provided by basic graphs of symmetric association schemes
and in particular by distance regular and strongly regular graphs
\cite{BanI84,BroCN89,PetR01}.
A serious attempt to establish a more strict approach to algebraic properties of
integral graphs is presented in \cite{SteAFV07}.
Clearly, rational circulants form an interesting particular case of regular integral graphs.
Investigation of these graphs usually is based on the amalgamation of techniques from
number theory, linear algebra and combinatorics.  Even a brief glance on such recent
contributions as \cite{So05,SaxSS07,AhmABS09, KloS10} shows a promising potential to
use for the same purposes also S-rings.

Let us now consider a very particular infinite series of rational
circulants $X_n=\cay(\Z_n,\Z_n^*)$, that is, the basic graph of the complete
rational circulant association scheme, containing edge $\{0,1\}$.
As in \cite{DejG95}, we will call such graphs {\em unitary} circulant graphs.
Different facets of interest to the unitary circulants may be traced from
\cite{Fuc05,KloS07,AkhBJJKKP09,RamV09,Dro10}.
A problem of description of $\aut(X_n)$ was posed in \cite{KloS07} and solved
in \cite{AkhBJJKKP09}.
Clearly, the reader will understand that the answer was in fact known for a
few decades in framework of the approach presented in this paper.
Similarly, one sets complete answer on the Problem 2 from \cite{KloS07}.

\subsection{Designed experiments: a bridge from and to statisticians.} 

I. Schur and R.~C.~Bose are now commonly regarded as the two most influential for-runners of
the theory of association schemes, a significant part of algebraic combinatorics,
see e.g. \cite{BanI84,KliRRT99,Bai04}.

A geometer by initial training, Bose (1901-1987) was in a sense recruited by P.~C.~Mahalanobis
to start from the scratch research in the area of statistics at a newly established
statistical laboratory at Calcutta (now the Indian Statistical Institute).
Fruitful influence of R.~A.~Fischer and F.~Levi (during 1938 - 1943 and later on)
turned out to become a great success not only for Bose himself, but also for all
growing new area of mathematics, see \cite{Bos82}.  As a result, within about two decades,
theory of association schemes was established by Bose et al.,
see \cite{BosN39,BosS52,BosM59,Bos63}
for most significant cornerstone contributions on this long way.
Being in a sense a mathematical bilingual, Bose was perfectly feeling in
the two areas which were created and developed via his very essential contributions:
design of statistical experiments and association schemes.

Unfortunately, over the theory of association schemes was recognized as an
independent area of mathematics, in particular after death of Bose, close links of
algebraic combinatorics to experimental statistics became less significant,
especially in the eyes of pure mathematicians.
Sadly this divergence still continues.
Nevertheless, mainly to the efforts of R.~A.~Bailey, a hope for the future reunion
is becoming during the last years more realistic.
The book \cite{Bai04} is the most serious messenger in this relation.
Being also bilingual (Bailey got initial deep training in classical group theory),
during last three decades she systematically promotes better understanding of
foundations of association schemes by statisticians. Referring to \cite{Bai04} for
more detail, we wish just to cite here such papers as \cite{Bai81,SpeB82,Bai85}
and especially \cite{Bai96}.

These contributions, became in turn, very significant for pure mathematics.  Indeed,
initial ideas of Nelder \cite{Nel65}, equivalent in a sense to the use of simple
reduction rules, in hands of Bailey et al. were transformed to the entire theory of
orthogonal partitions, group poset structures and crested products.
Note also that our striking example appears in \cite{Bai04} as Example~9.1 
in surprising clothes of designed experiment for bacteria search in a milk laboratory.

\subsection{Lattices and finite topological spaces.} 

For a square free number $n$ Gol'fand established in \cite{Gol85} bijection between
rational S-rings over $\Z_n$ and finite topologies on a $k$-element set,
here $n$ has exactly $k$ distinct prime factors.  This is a particular case of a
bijection between rational S-rings over $\Z_n$ and sublattices of $L(n)$ for arbitrary $n$.
Here we face another impact of diverse techniques from algebraic combinatorics,
general algebra, group theory, experimental designs, etc.
Such references as \cite{Sta71,Har74,Bel79,Duq86,Kod94,Ren94,Pfe04} provide a
possibility to make a brief glance of the top of this iceberg.

\subsection{Generalized wreath products.} 

The operation of wreath product has a long history, which goes back to such names as
A.~Cauchy, C.~ Jordan, E.~Netto and Gy.~P\'olya. E.~Specht was one of the first experts who considered it 
in a rigorous algebraic context, see \cite{Spe33}.
A new wave of interest and applications of wreath products was initiated
by L.~A.~Kalu\v{z}nin.  The Kalu\v{z}nin-Krasner Theorem (see \cite{KraK50-51}) is
nowadays commonly regarded as a classical result in the beginning course of group theory.
Less known is a calculus for iterated wreath product of cyclic groups, the outline of which
was created by Kalu\v{z}nin during the period 1941-45 (at the time
he was imprisoned in a nazi concentration camp), see \cite{Sus98}.
After the war the results, shaped mathematically, were reported on the Bourbaki seminar, and published in a series
of papers, see e.g. \cite{Kal48}. A few decades later on this calculus was revived, extended and
exploited in hands of L.~A.~Kalu\v{z}nin, V.~I.~Sushchaskii  and their disciples, cf. \cite{KalS73}.
The notation, used in current paper is inherited from the texts of Kalu\v{z}nin et al.

The generalized wreath product, the main tool in the reported project, was created
independently, more or less at the same time by two experts.
The approach of V.~Fe\v{\i}nberg (other spelling is Fejnberg) has purely combinatorial origins, first it was presented on the
IX All Union Algebraic Colloquium (Homel, 1968, see \cite{Fei68}). Details are given in a series
of papers \cite{Fej69a,Fej69b,Fej71,Fej73}. Fe\v{\i}nberg traces roots of his approach to the ideas of
Kalu\v{z}nin \cite{Kal51}. The book \cite{KalBF87} provides a helpful detailed source
for the wide scope of diverse ideas, related to different versions of wreath products,
its generalizations and applications. It seems that as an entity this stream of
investigations is overlooked by modern experts.

W.~Ch.~Holland submitted his influential paper \cite{Hol69} on January 11, 1968.
Though his interests are of a purely algebraic origin and the suggested operation is
less general (in comparison with one considered by Fe\v{\i}nberg), his ideas got much more
lucky fate. The paper \cite{Hol69} is noticed already in \cite{Wel76} and exploited in spirit of group posets
in \cite{Sil77} (both authors cite also \cite{Fej71}). It was Bailey who realized in
\cite{Bai81} that the approach of Holland is well suited for the description of
the automorphism groups of poset block structures.
With more detail all necessary main ideas may be detected from
\cite{BaiPRS83}, while \cite{PraRS85,Cam87} stress extra helpful information.
Our paper is strongly influenced by presentation in \cite{Bai04,Bai06}.

\subsection{Other products.} 

The crucial input in \cite{Bai06} is that the generalized wreath product of
permutation groups is considered in conjunction with the wreath product of association schemes,
and lines between the two concepts are investigated. The crested product is a particular case
of generalized wreath products, which may be alternatively explained in terms of iterated use
of crested product.  Note that the crested product for a particular case of S-rings was
considered in \cite{HirM01} under the name {\em star product}.
As we now are aware, the considered operations are enough in order to
classify rational S-rings over cyclic groups.

A more general product operation, the wedge product of association schemes, was recently
introduced and investigated in \cite{Muz09}. The term goes back to \cite{LeuM96,LeuM98}
who used it for recursive classification of S-rings over cyclic groups.
In a similar situation Evdokimov and Ponomarenko \cite{EvdP02}
are speaking about wreath product of S-rings
(unfortunately, their terminology does not coincide with traditional one).
Muzychuk also investigates the automorphism groups of his wedge product of association schemes.
For a particular case of S-rings over cyclic groups of prime-power order these groups
coincide with the subwreath  product in a sense of \cite{Kli94,KliP,Kov02,Kov}.
A few other operations over association schemes (semi-direct product and exponentiation)
are also of a definite interest, see references in \cite{Muz09},
though out of scope in this paper.

\subsection{More references.} 

It is a pleasure to admit that S-rings are proving their efficiency in algebraic graph
theory. As was mentioned, sometimes they may substitute the use of CFSG.
One more such example is provided by the classification of arc-transitive circulants.
This problem was solved for a particular case in \cite{XuBS04}, and in general
in \cite{Li05}. Both papers rely on a description of $2$-transitive groups
(a well known consequence of CSFG).
A particular case of such groups is considered in \cite{Jon02}.
In fact, the entire result in \cite{Li05} is a consequence of \cite{Muz94},
the proof runs in the same fashion as the one for Zibin's conjecture.

Below is a small sample of other situations when knowledge of S-ring theory turn out to be
quite helpful.

\begin{itemize}
\item Rational circulants, satisfying $A^{m} = d I + \lambda J$ \cite{Lam75,Ma84}.
\item Isomorphisms and automorphisms of circulants \cite{HuaM96,Ker09}.
\item Classification of distance regular circulants \cite{MikP03}.
\item Commuting decompositions of complete graphs \cite{AkbH07}.
\end{itemize}

For purely presentational purposes we also recall one more old example.
Arasu et al. posed in \cite{AraJMP94} a question about the existence of a Payley type
Cayley strongly regular graph $\Gamma$ which does not admit regular elementary abelian
subgroups of automorphisms. Such an example on $81$ points was presented in \cite{Kli94b}
as a simple exercise via the use of rational S-ring over group $\Z_{9}^2$,
it has automorphism group of order $1944$. An infinite series of similar examples,
using alternative techniques, was given \cite{Dav94}, automorphism groups were
not considered. Complete classification of such strongly regular graphs over $\Z_n^2$
with the aid of S-rings, is given in \cite{LeiM05} for $n=p^k$.
In our eyes the problem of classification of partial difference sets
(that is, Cayley strongly regular graphs) over groups $\Z_n^2, n \in \N$ is a nice
training task for innovative applications of S-rings and association schemes.

\subsection{Concluding remarks.} 

This project has been started in 1994 at the time of a visit of M.~Klin to Freiburg.
During years 1994-96 Klin was discussing with O.~H.~Kegel diverse aspects of the use of
S-rings and simple reduction rules. These discussions as well as ongoing 
numerous conversations with Muzychuk shaped the format of the project.
Starting from year 2003, Kov\'acs joined Klin, and by year 2006, in principle,
the full understanding of the automorphism groups of the rational circulants was achieved, and presented 
in \cite{KliK06}.
At that time we became familiar with \cite{BaiC05,Bai06}
and were convinced that the crested products is a necessary additional brick which allows
to create a clear and transparent vision of the entire subject.
Finally, a more ambitious lead was attacked; the authors were
striving to make presentation reasonably available to a wide mathematical audience.
Our goal is not only to solve a concrete problem but also to promote use of S-rings and to
stimulate interdisciplinary dialog between the experts from diverse areas, who for
many decades were working in a relative isolation, being not aware of the existence of
worlds ``parallel'' to their efforts.

\section*{Acknowledgements}

A visit of the second author at Ben-Gurion University of the Negev in November 2006 helped us to proceed with
the paper, the second author thanks Ben-Gurion University of the Negev for supporting his trip. 
The authors are much obliged to Otto Kegel and Misha Muzychuk for a long-standing fruitful cooperation. 
We thank Ilia Ponomarenko for helpful discussions and Andy Woldar for permanent stimulating interest to 
diverse facets of S-ring theory. Helpful remarks of Valery Liskovets are greatly appreciated.

\end{document}